\documentclass[12pt,leqno]{amsart}

\def \R {\Bbb R}
\def \cal {\mathcal}

\def \vs{\vspace*{0.1cm}}

\def \ds{\displaystyle}

\def\wt{\widetilde}
\def\i{\infty}
\def\D{\Delta}
\def\p{\partial}
\def\l{\lambda}
\def\ve{\epsilon}
\def\s{\sigma}
\def\a{\alpha}
\def\b{\beta}
\def\t{\theta}

\def\cal{\mathcal}

\def\l{\lambda}
\def\n{\nabla}
\def\O{\Omega}
\def\be1{{\begin{equation}}}
\def\ee1{{\end{equation}}}

\def\D{\Delta}
\def\n{\nabla}
\def\p{{\partial}}

\def\C{\mathbb C}

\def\l{\lambda}
\def\L{\Lambda}
\def\b{\beta}

\def\a{\alpha}

\def\e{\epsilon}
\def\d{\delta}
\def\S{\Sigma}

\def\part{\partial}
\def\R{{\mathbb R}}

\def\ba{\begin{array}}
\def\ea{\end{array}}

\newtheorem{them}{Theorem}[section]

\newtheorem{corollary}{Corollary}[section]

\newtheorem{pro}[them]{Proposition}

\numberwithin{equation}{section}

\newtheorem{lemma}{Lemma}[section]
\newtheorem{proposition}[lemma]{Proposition}
\newtheorem{theorem}[lemma]{Theorem}
\newtheorem{remark}[lemma]{Remark}

\title[]
{Analytic Aspects of the Toda System: II. Bubbling behavior and
existence of solutions}

\author[]{J\"{u}rgen Jost}
\address{Max Planck Institute for Mathematics in the Sciences, 04103 Leipzig, Germany}
\email{jost@mis.mpg.de}

\author[]{Chang-Shou Lin}
 \address{Department of Mathematics, National Chung-Cheng University,
 Minghsiung, Chia-Yi, Taiwan}
     \email{cslin@math.ccu.edu.tw}

\author[]{Guofang Wang}
 \address{Max Planck Institute for Mathematics in the Sciences, 04103 Leipzig, Germany}
     \email{gwang@mis.mpg.de}

\date{}

\begin{document}
\maketitle


\section{Introduction}
\setcounter{equation}{0}

In this paper, we continue to consider
the 2-dimensional (open) Toda system (Toda lattice)
for $SU(N+1)$
\begin{equation}\label{1.1}
- \D u_i = \sum_{j=1}^N a_{ij} e^{u_j},
\end{equation}
for $i=1,2,\cdots, N$, where $K=(a_{ij})_{N\times N}$ is the Cartan
matrix for $SU(N+1)$ given by
\[
\left(\begin{array}{rrrrrr}
2 & -1 & 0 & \cdots & \cdots & 0\cr
-1& 2 &-1& 0 & \cdots & 0\cr
0&-1&2&-1& \cdots & 0\cr
\cdots&\cdots&\cdots&\cdots&\cdots&\cdots\cr
0 &\cdots&\cdots& -1&2&-1\cr
0&\cdots&\cdots & 0 &-1&2\cr\end{array}\right).\nonumber
\]

As a very natural generalization of the Liouville
equation
\begin{equation}\label{1.2}
-\D u =2e^{u},
\end{equation} system (\ref{1.1})
is completely integrable, which
 is well-known in integrable systems theory.  The Liouville equation and the Toda
system arise in many physical models. In
Chern-Simons theories, the Liouville equation is closely  related
to Abelian models, while the Toda system is related
to non-Abelian models. See for instance
 the books \cite{dunnebook} and \cite{Yang}
and the references therein.

Though the Liouville equation has been extensively studied after Liouville
\cite{Liou}, the precise
behavior of convergence of its solutions has been  rather well
understood only in the last two decades, see for example \cite{BM, BLS,
ChenLi,
ChL, ChL2, DJLW, Li, LiS, NT}.
Such a delicate study leads to many
applications in the Abelian Chern-Simons theories and mean field equations.

To well understand the non-Abelian Chern-Simons model, we have to
study the analytic aspects of the Toda system. It is natural to
ask if we can generalize delicate analytic results for the
Liouville equation to the Toda system. However, the analysis of
the Toda system becomes more difficult, because the basic analytic
tool, the maximum principle, does not work. In  \cite{JW1}, we
established a Moser-Trudinger type inequality for the Toda system,
while a rough bubbling behavior of solutions to  system
(\ref{1.1}) was considered. See also \cite{LN}. Its bubbles
-entire solutions- were classified in \cite{JW2}. In this paper,
we will prove existence results for solutions to the Toda
system in various cases
by using methods developed in \cite{DJLW, NT, ChL, St, ST, Li}.
And we will give a much more
precise bubbling behavior of solutions, which will be very useful in our further
study of the Toda system.

Let $\S$ be a Riemann surface with
Gaussian curvature $K$. We consider the
following system
\begin{equation}
\label{1.3}\begin{array}{rcll} -\D u_i &=&\ds\vs \sum^N_{j=1}
\rho_ja_{ij}\left(\frac{h_je^{u_j}} {\int_\S h_je^{u_j}
}-1\right),&\quad \hbox{ in }
 \Sigma\quad 1\le i\le N,\\
\end{array}
\end{equation}
for the  coefficient matrix $A=(a_{ij})_{N\times N}$, the  Cartan
matrix of $SU(N+1)$ and $\rho=(\rho_1,\rho_2,\cdots, \rho_N)$ with
$\rho_i>0$ ($i=1,2,\cdots, N$) given  constants. Here $h_i:\Sigma \to \R$
is a given positive $C^1$ function for $i=1,2\cdots,N$. System
(\ref{1.3}) is the Euler-Lagrange system of the 
functional
\begin{equation}\label{1.4}
J_\rho(u)=\frac12\sum_{i,j=1}^N\int a^{ij}\n u_i\n
u_j+\sum_{j=1}^N\int \rho_ju_j-\sum_{j=1}^N \rho_j
\log\int_{\S}h_je^{u_j},
\end{equation}
in $H:=(H^1(\S))^N$ where $(a^{ij})$ is the inverse matrix of $A$.
In \cite{JW1}, we proved that $J_\rho$ has a lower bound  in $H$
if and only if $\rho_i\le 4\pi$ for any $i$, which is the Moser-Trudinger
inequality for the Toda system. From this inequality, we know that
if $\rho_i<4\pi$ ($\forall j$) then
$J_\rho$ iscoercive condition and hence $J_\rho$ has a minimizer,
which certainly satisfies system (\ref{1.3}).
 When one of the $\rho_i$'s equals  $4\pi$, the existence problem
 becomes subtler. For simplicity of notation, we consider only the case $N=2$.

  Our first result in this case is
\begin{theorem}\label{thm1.1} Let $\S$ be a Riemann surface with
Gaussian curvature $K$ and $N=2$. And let $\rho_1=4\pi$ and
$\rho_2\in (0, 4\pi)$. Suppose that
\begin{equation}\label{1.5} \D \log h_1(x)+(8\pi-\rho_2)-2K(x)>0\ \ \mbox{for}\ \ x\in\S.
\end{equation}
 Then $J_\rho$ has a minimizer
$u=(u_1,u_2)$ satisfying
\begin{equation}\label{1.6}
\begin{array}{rcl}
-\D u_1&=&\ds\vs 2\rho_1\left(\frac {h_1e^{u_1}}{\int_\S
h_1e^{u_1}}-1\right)
-\rho_2\left(\frac {h_2e^{u_2}}{\int_\O h_2e^{u_2}}-1\right)\\
-\D u_2&=&\ds 2\rho_2\left(\frac {h_2e^{u_2}}{\int_\S
h_2e^{u_2}}-1\right) -\rho_1\left(\frac {h_1e^{u_1}}{\int_\S
h_1e^{u_1}}-1\right),
\end{array}\end{equation}
for $\rho_1=4\pi$ and $\rho_2\in (0, 4\pi)$.
\end{theorem}

Theorem \ref{thm1.1} is obtained by using a result in \cite{ChL}, which
is a refined result of \cite{DJLW} and \cite{NT}. For $\rho_1=\rho_2=4\pi$,
we also have a sufficient condition under which (\ref{1.6}) has a solution.
See Theorem \ref{thm6} below.

For the general case, we have the following compactness of the solution space.

\begin{theorem}\label{thm1.2} For any compact set $\L_1\times \L_2 \subset \R_+\times
\R_+$, if there are two positive integers $m_1$ and $m_2$ such that $\L_i \subset
(4\pi m_i, 4\pi (m_i+1))$ for $i=1,2$, then the
solution space of (\ref{1.6}) for $(\rho_1, \rho_2)\in \L_1\times \L_2$
is compact.
\end{theorem}
The proof of Theorem \ref{thm1.2} follows from the study of the convergence
of the sequence of solutions, which was initiated in \cite{JW1} for the Toda
system.
Let $x$ be an  blow-up point of the sequence
$u^k$, i.e., there exists a sequence $x_k\to x$ such that
$\max\{u_1^k(x_k), u_2^k(x_k)\}\to \infty$ as $k\to \infty$.
Define
\[\s_i=\lim_{r\to 0}\lim_{k\to \infty}\int_{B_{r}}e^{u^k_i}.\]
We call $(\s_1, \s_2)$ a blow-up value at $x$. First it is easy to check that
$\s_1+\s_2>0$. A local Pohozaev argument gives us a relation between
$\s_1$ and $\s_2$.
$$\s_1^2 + \s_2^2 -\s_1 \s_2 = 4\pi (\s_1+\s_2).$$
See \cite{LN} and \cite{JW1} for the proof.
In Proposition \ref{pro} below,
we show that $(\rho_1, \rho_2)$ can  only be one of
$(4\pi, 0)$, $(0,4\pi)$, $(4\pi, 8\pi)$, $(8\pi, 4\pi)$ and
$(8\pi, 8\pi)$. It is clear that Theorem \ref{thm1.2} is a direct consequence of
Proposition \ref{pro}. In the first two cases, one of the $u_i^k$ does not blow-up,
another bubbles more or less like  solutions of the Liouville equation,
which has been extensively studied in the last decade.
The last case contains a typical blow-up phenomenon of solutions of the Toda system.
It is one of our main results to obtain a precise description of the bubbling behavior for this case.
To describe it, we assume
that there is a sequence of solutions $u^k=(u^k_1,u^k_2)$ of
\begin{equation}\label{1.7}
\left\{
\begin{array}{rcl}
-\D u^k_1 &=& \vs\ds 2 h_1^ke^{u^k_1}-h^k_2e^{u^k_2}\\
-\D u^k_2 &=& \ds 2 h_2^ke^{u^k_2}-h^k_1e^{u^k_1}
\end{array} \quad \hbox{ in } B_2.\right. \end{equation}
Here $B_r$ denotes a disk of radius $r$  and center $0$ and $h^k_i$ converges in
$C^1$ to a positive $C^1$ function $h_i$ for $i=1,2$. Assume
without loss of generality that $h_1(0)=h_2(0)=1$. Suppose that
$u^k$ bubbles off, i.e., $\max_{x\in B_2} \{u^k_1, u^k_2\}\to
\infty$ as $k\to \infty$. More precisely we assume that
\begin{enumerate}
\item[(1)] 0 is the only blow-up point of $u^k$. \item[(2)]
$\max_{\p B_2} u_i^k - \min_{\p B_2} u_i^k \leq c$ for $i=1,2$.
\item[(3)] $ \int_{B_2} e^{u_i^k} dx \leq c$ for $i=1,2$ and any
$k$.
\end{enumerate}
Assume that $\l^k = \l_1^k := \max_{B_2} u_1^k \geq \max_{B_2}
u_2^k =: \l_2^k$. Let $x_k$ be the maximum point of $u_1^k$. Set
$v^k=(v_1^k, v_2^k)$ by
$$v_i^k(x) = u_i(\ve_k x+x^k)-\l^k\  \ \mbox{for}\ \ i=1,2,$$
where $\ve_k=e^{-\frac 1 2 \l_k}$. Clearly $v^k$ satisfies in
$\{ x\in \R^2 \mid \ve_k x+x^k \in B_2\}$
$$\left\{ \begin{array}{l}
\ds -\D v_1^k = 2h^k_1(\e_kx+x^k) e^{v_1^k} -  h^k_2(\e_kx+x^k)e^{v_2^k},\\
\ds -\D v_1^k = 2 h^k_2(\e_kx+x^k)e^{v_2^k} -  h^k_1(\e_kx+x^k)e^{v_1^k},\\
\end{array}\right.$$
and $v_i^k\leq 0, v_1^k(0)=0$. Since we only consider the case
that the bubble is a solution of the Toda system, we may further assume that
\begin{enumerate}
\item[(4)] $v_2^k(0)$ is bounded from below.
\end{enumerate}
Then, there exists a solution $v^0=(v_1^0, v_2^0)$ of the Toda
system
$$\left\{ \begin{array}{l}
-\D v_1^0 = 2e^{v_1^0} - e^{v_2^0},\\
-\D v_2^0 = 2e^{v_2^0} - e^{v_1^0},\\
\end{array}\right.$$
such that $v_i^k-v_i^0$ converges to zero in $C_{loc}^2(\R^2)$.

{From} the classification result for entire solutions of the Toda system
\cite{JW2}, which is a generalization of the classification result
for the entire solutions of the Liouville equation obtained by
Chen and Li \cite{ChenLi}, $v^0$ is obtained from a rational curve
from $S^2$ to ${\mathbb C} P^2$. In particular, we have
\begin{equation}{\label{1.8}}
\int_{\R^2} e^{v_1^0} = \int_{\R^2} e^{v_2^0} =8\pi.
\end{equation}
 Now we state our main theorem.

\begin{theorem}\label{thm1.3} Let $u^k=(u^k_1,u^k_2)$ be a
sequence of solutions to (\ref{1.7}). Suppose that
 $h^k_i$  converges in $C^1$ to a positive function $h_i$ with $h_i(0)=1$ for $i=1,2$
and  that $(1)$-$(4)$ hold. Then there exist two constants $r_0>0$
and $c>0$ independent of $k$, such that
\begin{equation}{\label{1.9}}
|u_i^k(x)-\l^k-v_i^0 (\ve_k^{-1}(x-x^k))|<c\ \ \mbox{in}\ \
B_{r_0}
\end{equation}
for $i=1,2$.
\end{theorem}

Theorem \ref{thm1.3} gives a precise asymptotic behavior of a
blow-up sequence of solutions. When $N=1,$ Theorem \ref{thm1.3} was
proved by Y.Y Li \cite{Li} by using the method of moving planes
 and a symmetry of the
entire solutions of the Liouville equation. As mentioned above,
the method of moving planes, which  is based on  the maximum
principle, does not work for the Toda system. And also, the
symmetry used in \cite{Li} is not available for the entire
solutions of the Toda system in general. Thus the method in
\cite{Li} could not provide a proof for the Toda system. Instead,
 here, together
with the geometry related  to the Toda system, we will employ a more
delicate analysis. This is a combination of methods given in
\cite{CLT} and \cite{JW2}.

Our last result is an existence result for a supercritical case with
respect to the Moser-Trudinger inequality established in \cite{JW2}.

\begin{theorem}\label{thm1.4} Let $\Sigma$ be a compact Riemann surface of genus greater than $0$.
Then for any $\rho=(\rho_1, \rho_2)$ with
$\rho_i\in (0,4\pi)\cup (4\pi, 8\pi)$
system (\ref{1.6}) has a solution.\end{theorem}

The method to prove Theorem \ref{thm1.4} follows \cite{St}. See also
\cite{St-book, ST, DJLW} and also \cite{LN}. One can also prove  the Theorem
by computing a topological degree of system (\ref{1.6}). In fact,
the results -especially Theorem \ref{thm1.3}- established in this paper
will be used in computing the  topological degree for the Toda system. We will pursue this subject in
forthcoming papers.

In Section 2, we first recall  basic facts about the convergence
of  solutions to the Toda system and   we then show that the
possible blow-up values of the Toda system are isolated.
 We recall the geometric
interpretations of solutions to the Toda system in Section 3. Such
a solution corresponds to a flat connection and its singularity to
the holonomy of its corresponding flat connection. We will give a
more precise bubbling behavior, Theorem \ref{thm1.3}, in Section 4.
The proof is a fine combination of arguments presented in
\cite{CLT} and \cite{JW2}. In Section 5, we prove Theorem \ref{thm1.1}
by applying a result in \cite{ChL}.
In Section 6, we give an existence
result for (\ref{1.1}) with the Dirichlet boundary for a supercritical
case, based on a method of Struwe \cite{St} and the
Moser-Trudinger inequality established in \cite{JW1}.  We give an asymptotic
behavior of singularities of solutions in Appendix.

\noindent{\it Acknowledgement.} A part of the paper was carried out
while the third author was visiting the Center of Theoretical
Science in Taiwan (CTS). He would like to thank the CTS for warm hospitality.

\section{Bubbling behaviors}

The convergence of solutions of Toda system was studied in
\cite{JW1}. Here for the convenience of the reader, we first
recall some basic facts.

\begin{proposition}\label{P2.1}
Let $\Omega$ be a bounded smooth domain in $\R^2$
and $u^k=(u^k_1,u^k_2)$ be a sequence of solutions of the
following system
\begin{equation}\label{eq100}
\left\{\begin{array}{rcll}
-\D u^k_1& =& 2h_1^k e^{u^k_1}-h_2^k e^{u^k_2}, &\text{ on } \O,\\
-\D u^k_2& =& 2h_2^k e^{u^k_2}-h_1^k e^{u^k_1}, &\text{ on } \O,
\end{array}\right.\end{equation}
with
\begin{equation}\label{3a}
\int_\O e^{u^k_1}<C,\ \text {and}\  \int_\O
e^{u^k_2}<C.
\end{equation}
Set
\begin{equation}\label{3b}
 S_j=\{x\in \S| \exists \text{ a sequence $y^k\to x$ such that }
u^k_j(y^k) \to +\infty\}.\end{equation} Then, one of the
following possibilities happens: (after taking subsequences)
\begin{itemize}
\item [(1)] $u_i^k$ is bounded in $L^\infty_{loc}(\O)\times
L^\infty_{loc}(\O)$. \item [(2)] For some $j\in\{1,2\}$, $u^k_i$
in $L^\infty_{loc}(\O)$, but $u^k_j\to -\infty$ uniformly on any
compact subset of $\O$ for $j\not = i$. \item[(3)] For some $i\in
\{1,2\}$, $S_i\not = \emptyset$, but $S_j = \emptyset$, for
$j\not = i$. In this case, $u^k_i\to -\infty$ on any compact
subset of $\O{\backslash} S_i$, and either, $u^k_j$ is bounded in
$L^\infty_{loc}(\O)$, or $u^k_j \to -\infty$ on any compact subset
of $\O$. \item[(4)] $S_1 \not =\emptyset$ and $S_2 \not
=\emptyset$. Moreover, $u^k_j$ is  either bounded or $\to -\infty$
on any compact subset of $\O{\backslash} (S_1\cup\S_2)$ for
$j=1,2$.
\end{itemize}
\end{proposition}

\begin{proof} The proof is given in \cite{JW1}. The idea follows
closely \cite{BM}.\end{proof}

\begin{remark} \label{rem3.2} One can prove that when $S_1\neq \emptyset$, $S_2\neq \emptyset$ and
 $S_1\backslash S_2\neq \emptyset,$  $u^k_1\to -\i$ uniformly in any compact set of
$\O\backslash\{S_1\cup S_2)$. See \cite{BM, JW1}. However, the proof of
Proposition \ref{pro} below implies that when
$S_1\neq \emptyset$ and $S_2\neq \emptyset$, both $u_1^k\to -\i$
uniformly in any compact set of $\O\backslash\{S_1\cup S_2)$. This is an improvement of (4).
\end{remark}

In this paper, we do not distinguish convergence and subconvergence.
We may assume that there exist two nonnegative bounded measures
$\mu_1$ and $\mu_2$ such that
\[e^{u^k_i}\psi \to \int \psi d\mu_i \text{ as } k\to \infty,\]
for every smooth function $\psi$ with support in $\O$ and $i=1,2$.
A point $x\in \O$ is called a {\it $\gamma$-regular point} with
respect to $\mu_j$ if there is a function $\psi\in C_c(\O)$ ,
$0\le \psi \le 1$, with $\psi=1$ in a neighborhood of $x$ such
that
\[\int_\O \psi d\mu_j <\gamma.\]
We define
\[\O_j(\gamma)=
\{x\in\O\,|\,x \text{ is not a }\text{$\gamma-$regular point with
respect to } \mu_j\}.
\]
One can show  $\O_1(\gamma)$ and $\O_2(\gamma)$ are finite. And
$\O_j(\gamma)$ is independent of $\gamma$ for small $\gamma<
2\pi$, see \cite{JW1}. Furthermore, we have
\begin{equation}\label{eq3c}
S_i=\O_i(\gamma),\quad \hbox{ for }\gamma<2\pi.\end{equation}

Proposition \ref{P2.1} implies that  the blow-up points are
isolated. Let $0$ be an isolated  blow-up point of the sequence
$u^k$, i.e., there exists a sequence $x_k\to 0$ such that
$\max\{u_1^k(x_k), u_2^k(x_k)\}\to \infty$ as $k\to \infty$.
Define
\[\s_i=\lim_{r\to 0}\lim_{k\to \infty}\int_{B_{r}}e^{u^k_i}.\]
We call $(\s_1, \s_2)$ a blow-up value. If one of $\s_1$ and
$\s_2$ is zero, then another is $4\pi$.
This can also be obtained from the following local Pohozaev
identity.

\begin{lemma}\label{lem3.1} We have
$$\s_1^2 + \s_2^2 -\s_1 \s_2 = 4\pi (\s_1+\s_2)$$\end{lemma}
\begin{proof} See the proof in \cite{JW1} and \cite{LN}.
\end{proof}

Hence we assume that $\s_i>0$ for $i=1,2$. From (\ref{eq3c}) we
know that $\s_i\ge 2\pi$ for $i=1,2$. This can be improved in the
following
\begin{lemma}\label{lem3.0} We have $\s_i\ge 4\pi$ for $i=1,2$.
\end{lemma}
\begin{proof} A proof of this Lemma is given in \cite{LN}.
Here for convenience of the reader, we give a proof. Assume by
contradiction that $\s_1<4\pi$. From Lemma \ref{lem3.1},
we have $\s_2 <8\pi$. Let $x_k$ be the maximum point of
$u^k_1$ and $\l_k=\max u^k_i$ the maximum value of $u^k_1$. Since
$\s_1>0$, we know that $\l_k\to \i$ as $k\to \i$. Consider
\[\tilde u^k_i(x)=u^k_i(e^{-\l_k/2}x+x_k)-\l_k,\]
in $\O_k:=\{x\in\R^2\,|\, e^{-\l_k/2}x+x_k\in B_{r_0}\}$ for a
fixed $r_0>0$. It is clear that $(\tilde u^k_1, \tilde u^k_2)$
satisfies (\ref{eq100}) in $\O_k$ and  $\O_k$ converges to $\R^2$
as $k\to \infty$. Now we
have three possibilities:
\begin{itemize}
\item[(i)] $\tilde u^k_2$ blows up. Namely  there is a point
$y\in\R^2$ and a sequence $y_k$ such that $y_k\to y$ and $\tilde
u^k_2(y_k)\to +\i$ as $k\to\i$.
\item[(ii)] $\tilde u^k_2$
uniformly converges to $-\i$ in any compact set of $\R^2$.
\item[(iii)] $\tilde u^k_2$ is bounded from above and there is a point
$x$ such that $\tilde u^k_2(x)\ge -C$ for  some constant $C>0$
independent of $k$.\end{itemize}

We first consider case (i). We define $\tilde S_2$ to be the blowup set of $\tilde u_2^k$.
It is clear that $\wt S_2$ is finite. In fact,  one can show that the
number of $\wt S_2$ is $1$.


Since $\tilde u_1^k(x)\leq 0$, by using the equation for $\tilde u_2^k$ alone, we can apply a
result in [13] to conclude that for each blow-up point $y\in \tilde S_2$, the local mass around
$y$ would converge to $4\pi$. Thus
$$4\pi\#\tilde S_2\leq \int_{\O_k} e^{\tilde u_2^k} dx \leq \s_2 \leq 8\pi.$$
Hence $\# \tilde S_2=1$.

Suppose $\tilde S_2=\phi$. Then both $\tilde u_1^k$ and $\tilde u_2^k$ converge to $v_1$ and $v_2$, which
are a solution of the Toda system (\ref{1.1}) with $N=2$. But the mass formulas (\ref{1.8})
implies $\s_2\geq 8\pi$, which is a contradiction to the assumption $\s_2<8\pi$.

Now suppose $\tilde S_2=\{p\}$. Then it is easy to see $\tilde u_1^k$ converges to $u$ in
$\R^2\backslash \{p\}$, and $u$ satisfies
$$-\D u = 2e^u - 4\pi \d_p.$$

Let $v=u-2\log|x-p|$. It is clear that $v$ satisfies
\[-\D v =2|x-p|^2e^v.\]
A result given in \cite{ww} (see also \cite{CLT}) gives that
\[2\int_{\R^2}e^u=\int_{\R^2}|x|^2e^v=16\pi,\]
which implies that $\s_1\ge 8\pi$. This contradicts $\s_1<4\pi$.

 In case (ii), one can show that
$\tilde u_1^k$ converges in $C^2_{loc}(\R^2)$ to a $u$, which is
an entire solution of (\ref{1.2}) with finite energy, the
classification result of Chen-Li in \cite{ChenLi} gives us
\[ \lim \int_{B_{r_0}} e^{u^k_1}=\lim \int_{\O_k}e^{\tilde u^k_1}\ge
4\pi,\] which contradicts $\s_1< 4\pi$. In case (iii), we can
show that $(\tilde u^k_1, \tilde u^k_2)$ converges to an entire
solution of the Toda system (\ref{1.1}) with finite energy. Now the
classification result for the Toda system \cite{JW1} gives us that
\[ \lim \int_{B_{r_0}} e^{u^k_1}=\lim \int_{\O_k}e^{\tilde u^k_1}\ge
8\pi,\]
again, a contradiction.
 \end{proof}

 {From} Lemma \ref{lem3.1}, the dimension of the set of possible blow-up values
is less than or equal to one. We will  show  that the possible values of
$(\s_1,\s_2)$ are in fact isolated.  It might be not difficult to see that
$(4\pi, 0)$ and $(0, 4\pi)$ are possible values. And $(8\pi,
8\pi)$ is also a possible blow-up value. On $\R^2$, the solution
space of the Toda system is noncompact and the blow-up value is
$(8\pi, 8\pi)$. Now we have

\begin{proposition}\label{pro}
The blow-up value of the Toda system (\ref{eq100}) can only be one
of $(4\pi, 0)$, $(0,4\pi)$, $(4\pi, 8\pi)$, $(8\pi, 4\pi)$ and
$(8\pi, 8\pi)$.\end{proposition}

\begin{proof} We only need to exclude the case  that one of $\s_1$
and $\s_2$ is greater that $8\pi$. Assume by contradiction that
$\s_2>8\pi $. In view of Lemma \ref{lem3.1}, we have $\s_2<12\pi$ and $4\pi<\s_1<8\pi$.

Choose $r_k$ such that
\begin{equation}{\label{eq:2-5}}
\int_{B_{r_k}}e^{u^k_2}=8\pi.
\end{equation}
Since $\s_2>8\pi$, it is easy to check that $r_k \to 0$ as $k\to
\infty$. Consider
\[\tilde u^k_i(x)=u^k_i(2r_k x)-2\log 2r_k,\]
in $B_{r_k^{-1}}$. Let $\widetilde S_i$ be the blow-up set of
$\tilde u^k_i$.


If $\widetilde S_2=\emptyset$ and $\widetilde S_1\neq \emptyset$, then $\tilde
u_2^k$ converges in $C^2_{loc}(\R^2\backslash \widetilde S_1)$ to
a solution of
\begin{equation}\label{eq3} -\D v_0=2e^{v_0}-4\pi
\sum_{p\in\widetilde S_1} \d_p.\end{equation} It is clear that the
number of $\widetilde S_1$ is $1$. Otherwise, $\s_1\ge 8\pi$,
which is a contradiction to our assumption. Hence as above we have
\[ \int_{\R^2}e^{v_0} =8\pi.\]
However by (\ref{eq:2-5}),  $8\pi=\int_{\R^2}e^{v_0}> \int_{B_2}e^{v_0} \ge \lim
\int_{B_1} e^{\tilde u^k_2}=8\pi$, a contradiction. Thus $\tilde S_1=\phi$.
If $\widetilde S_2=\widetilde S_1=\emptyset$, then as discussed in the proof of Lemma \ref{lem3.0},
either $\tilde u^k_2$ converges to an entire solution $\tilde v_0$
of the Liouville equation and $\tilde u^k_1$ converges uniformly
to $-\i$ in any compact domain of $\R^2$, or
 $(\tilde u^k_1, \tilde u^k_2)$ converges to an entire solution of
 the Toda system. The latter case is clearly a contradiction to $\s_1<8\pi$
 and the former  leads to
 \[4\pi=\int_{\R^2} e^{\tilde v_0}>\int _{B_2} e^{\tilde v_0}\ge
 \lim_{k\to\i} \int_{B_{r_k^{-1}}} e^{\tilde u^k_2}=8\pi,\]
 a contradiction again. Hence $\widetilde S_2 \neq \emptyset$.
By rescaling a factor of close to 1 if necessary,
we may assume that $\{x\in \R^2\, |\,
|x|=1\}\cap (\widetilde S_1\cup \widetilde S_2)=\emptyset$,
\begin{equation}\label{eq-2a}
\tilde u^k_i(x)\to -\infty, \quad \hbox{ locally in $\R^2\backslash (\tilde S_1\cup \tilde S_2)$
as } k\to \infty,\end{equation}
for $i=1,2$ and
\begin{equation}\label{eq200} \int_{\{|x|\le 1\}}e^{u^k_2} \to
8\pi\quad \hbox{ as } k\to \infty.
\end{equation}
Here (\ref{eq-2a}) holds by Remark \ref{rem3.2}.
 Now we use the Kelvin transformation to define
\[v^k_i(x)=\tilde u^k_i\left(\frac x{|x|^2}\right)
-4\log|x|.\] Since the Toda system is conformally invariant,
$v^k=(v^k_1, v^k_2)$ satisfies the Toda system in $\R^2\backslash
B_{r_k}$.
We consider $v^k$ on $\Omega_k:= B_1\backslash B_{r_k}$. The
boundary of $\Omega_k$  consists of two components
\[\partial_1\O_k=\{|x|=1\}\quad \quad \hbox{ and } \quad \partial_2\O_k=\{|x|=r_k\}.\]
On $\partial_1\O_k$, $v^k_i\to -\infty$ as $k\to \i$. Let
$\mu_k=\max_{\O_k}v^k_2$. Since $\s_2-8\pi>0$, it is clear that
$\mu_k\to \i$ as $k\to\i$. Let $y_k\in \overline \O_k$ such that
$v^k_2(y_k)=\mu_k$. We claim that
\begin{equation}\label{eq3f}\frac{|y_k|}{r_k} \to \i \quad \hbox{ as }
k\to\i.\end{equation}

Suppose that $|y_k|/r_k\le d$ uniformly for some constant $d>1$.
We consider $\tilde u^k$ in $B_{r_k^{-1}}$. We first have by integrating
by parts, for any $r\ge 1$,
 \begin{equation}\label{eq3h}
 -\frac{d}{dr}\bar {\tilde u}^k_2(r)r =\frac
1{2\pi}(2\int_{B_r} \tilde h_1^k e^{\tilde u^k_2}-\int_{B_r} \tilde h_1^k e^{\tilde u^k_1})\ge
\frac 1{2\pi}(16\pi-\s_1)>4(1+\e_o)\end{equation}
 for some
constant $\e_0>0$. Here $\bar{\tilde u}^k_i(r)=\frac 1{2\pi
r}\int_{|x|=r}\tilde u^k_i d\s$ is the average of $\tilde u^k_i$
over $\{|x|=r\}.$

Let $f(r)=\bar{\tilde u}_2^k(r)+4\log r$. From above we have $f'(r) \le
-\e_0/r$. Hence, we have for $r\ge 1$ \begin{equation}
\label{eq3d} \bar {\tilde u}^k_2(r)+4\log r \le f(1)= \bar
{\tilde u}^k_2(1)=:c_k.\end{equation}
 {From} above, we know $c_k\to -\i$ as
$k\to \i$.
 Now by the definition of $y_k$, we have
for any $|x|\ge 1$ \begin{equation}\label{eq3e} \ba {rcl} \ds\vs
{\tilde u}^k_2(x)+4\log|x|&=&\ds
v^k_2(\frac x{|x|^2})\le  v^k_2({y_k})\\
&=&\ds {\tilde
u}^k_2(\frac{y_k}{|y_k|^2})-4\log|y_k|.\ea\end{equation}
Recall that  there is a constant $c_d>0$, independent of $k$, such
that \begin{equation}\label{eq3j}
|u^k_2(x)-u^k_2(y)|<c_d,\end{equation}
for $|x|\ge d^{-1},\,\,|y|\ge d^{-1}$. Thus
\begin{equation}{\label{eq:2-13}}
|\tilde u_2^k(x)- \tilde u_2^k(y)|\leq c_d
\end{equation}
for $|x|\geq r_k^{-1} d^{-1}$ and $|y|\geq r_k^{-1} d^{-1}$. (\ref{eq:2-13}),
together with $\frac {r_k}{|y_k|}\ge d^{-1}$, implies that
\begin{equation}\label{eq3j1}
|\bar {\tilde u}^k_2(\frac {y_k}{|y_k|^2})- {\tilde u}^k_2(\frac
{y_k}{|y_k|^2})|<c_d.\end{equation}
(\ref{eq3j1}), together with (\ref{eq3d}) and
(\ref{eq3e}), implies that for $|x|\ge 1$,
\[\ba{rcl}
\ds\vs {\tilde u}^k_2(x)+4\log|x| &\le & \ds \tilde
u^k_2(\frac{y_k}{|y_k|^2})-4\log|y_k|\\
&\le & \ds\vs \bar{\tilde u}^k_2(\frac{y_k}{|y_k|^2})-4\log|y_k| +c_d\\
&\le& \ds  c_k+c_d \to -\infty \ \text{as}\ k\to +\infty.\ea\]
Therefore, we have
\[\int_{1\leq |x|\leq \frac{1}{r_k}}e^{\tilde u_2^k}\le e^{c_k+c_d}
\int_{|x|\ge 1}|x|^{-4} =O(1)e^{c_k}\to 0\]
as $k\to\i$, which implies that $\s_2=8\pi$, a contradiction. Thus, (\ref{eq3f}) is proved.


Consider a new rescaled $\tilde
v^k=(\tilde v^k_1, \tilde v^k_2)$ defined by
\[\tilde v^k_i(x)=v^k_i(e^{-\mu_k/2}x+y_k)-\mu_k\]
in
\[\widetilde\O_k=\{x\in\R^2\,|\, e^{-\mu_k/2}x+y_k\in \O_k\}.\]
{From} (\ref{eq3f}), we have two possibilities: \begin{itemize}
\item [(i)] $|y_k|e^{\mu_k/2}\to \i$ as $k\to\i$. \item[(ii)]
$|y_k|e^{\mu_k/2}$ is uniformly bounded. \end{itemize}For case
(i), we know that
 $\widetilde \O_k \to \R^2$. Arguing as in the proof of Lemma \ref{lem3.0},
 we can show that $\tilde
v^k_2$ converges in a suitable topology to an entire solution to
an equation like (\ref{eq3}) or $(\tilde v^k_1,\tilde v^k_2)$ to
an entire solution of the Toda system. In both cases, we have
\[\lim \int_{\O_k}e^{\tilde u^k_2}\ge 4\pi,\]
which, it turn, implies that $\s_2\ge 12\pi$, a contradiction. For
case (ii), we assume that $y_ke^{\mu_k/2}\to p\in \R^2$. In view
of (\ref{eq3f}), we have that $\wt \O^k$ converges to
$\R^2\backslash\{p\}$. In this case, one can show as above that
$\tilde v^k_2$ converges to a function $\tilde v_0$ in
$C^2_{loc}(\R^2\backslash(\{p\}\cup\hat S_1))$, where $\hat S_1$
is the set of blow-up points of $\tilde v^k_1$.
One can check that $\tilde v_0$ satisfies
\begin{equation}\label{eq3g}
-\D \tilde v_0=2e^{\tilde v_0}-\a 4\pi \d_p-4\pi \sum_{x\in\hat
S_1} \d_x.\end{equation}
Since $\tilde v^k_2\le 0$, we have
$\a\ge 0$. It is clear that the number of $\hat S_1$ is less than
$2$ because $\s_1<8\pi$. If $\hat S_1=\emptyset$,  together with $\a\geq 0$,
then we have $\int_{\R^2}e^{\tilde v_0}\ge 4\pi$, which implies that $\s_2\ge 12\pi$, a
contradiction. Now assume that $\hat S_1=\{q\}$. We will show  in
the appendix  that $\int_{\R^2}e^{\tilde v_0}>4\pi$, a
contradiction again. This completes the proof of the Proposition.
\end{proof}

\section{Geometric interpretations}
Before we give a proof of the precise asymptotic behavior of solutions
of the Toda system,  we would like to recall the geometric interpretation of
solutions of the Toda system (\ref{1.1}).

Let $\Omega$ be a simply connected domain and $u=(u_1,u_2, \cdots,
u_N)$  a solution of (\ref{1.1}) on $\Omega$. Define $\tilde
w_0,\tilde  w_1,\tilde w_2,\cdots, \tilde w_N$ by
\begin{equation}
\label{2.1} u_i=2\tilde w_i-2\tilde w_{i-1} \quad\text{ for } i\in
I \text { and } \sum_{i=0}^N\tilde w_i=0.\end{equation} It is
clear that
\[\tilde w_0=-\frac 1{2(N+1)}\sum_{i=1}^N (N-i+1)u_i.\]
Now set $w_i=\tilde w_i-i\log 4$ for $i=1,2,\cdots,N$. It is easy
to check that $w_0, w_1,\cdots, w_N$ satisfies
\begin{equation}\label{2.3}\left\{ \begin{array}{lll}
2(w_0)_{z\bar z}& =& e^{2(w_1-w_0)}\\
2(w_1)_{z\bar z}& =&   -e^{2(w_1-w_0)}+  e^{2(w_2-w_1)}\\
 \cdots  & \cdots&\cdots\\
2(w_N)_{z\bar z} &=& -e^{2(w_N-w_{N-1})}.\\
\end{array} \right.\end{equation}

(\ref{2.3}) is equivalent to an integrability condition
\begin{equation}\label{integ}
{\cal U}_{\bar z}-{\cal V}_{ z}= [{\cal U}, {\cal V}]
\end{equation}
of the following two equations
\begin{equation}\label{2.4}
\phi^{-1}\cdot\phi_z= {\cal U}\end{equation} and
\begin{equation}\label{2.5}
\phi^{-1}\cdot\phi_{\bar z} ={\cal V,}\end{equation} where
\[{\cal U}=
\left(\begin{matrix} (w_0)_{z} & &&\\
& (w_1)_z && \\
&& \cdots & \\
&&& (w_N)_z \\
\end{matrix}\right)+
\left(\begin{matrix} 0 & e^{w_1-w_0} &&\\
 & 0  && \\
 && \cdots &e^{w_N-w_{N-1}} \\
&&& 0 \\
\end{matrix}\right)
\]
and
\[{\cal V} =
-\left(\begin{matrix} (w_0)_{\bar z} & &&\\
& (w_1)_{\bar z} && \\
&& \cdots & \\
&&& (w_N)_{\bar z} \\
\end{matrix}\right)-\left(\begin{matrix} 0 & &&\\
e^{w_1-w_0} & 0  && \\
 && \cdots & \\
&&e^{w_N-w_{N-1}}& 0 \\
\end{matrix}\right)
\]

Set

\begin{equation}\label{eq3.5a}
\alpha=-({\cal U}dz+{\cal V}d \bar z).\end{equation} Then,
$\alpha$ is a one-form valued in $su(N+1)$. With the help of the
Frobenius Theorem, we obtain a map $\phi:\O\to SU(N+1)$ such that
\[\a=\phi^{-1}\cdot d\phi.\]
As a connection on the trivial bundle $\Omega\times \C^{N+1}\to
\Omega,$
 $\a$ (or $d+\a$) is  flat, i.e.,
$\a$ satisfies the Maurer-Cartan equation
\[ d\a+\frac 12[\a,\a]=0.\]

When $\Omega$ is not simply connected, we cannot apply the
Frobenius theorem directly to obtain $\phi$. Let
$\Omega=B^*=B\backslash\{0\}$ be the punctured disk. We introduce
the {\it holonomy} (see \cite{SS}) of an $SU(N+1)$ connection
$\a$ on the bundle $\Omega\times \C^{N+1}\to \Omega$ along  around
$0$. Let $(r,\theta)$ be the polar coordinates. Decompose $\a$ as
$\a=\a_rdr+\a_{\theta} d\theta$. $\a_r$ and $\a_\theta$ are
$su(N+1)$-valued. For any given $r\in (0,1)$, the following
initial value problem,
\[ \frac {d\phi_r}{d\theta}+\a_{\theta}\phi_r =0, \quad \phi_r(0)=Id,\]
has a unique solution $\phi_r(\theta)\in  SU(N+1)$. Here $Id$ is
the identity matrix. If the connection $d+\a$ has no singularity
at $0$, the conjugacy class of $\phi_r(2\pi)$ is trivial. If $\a$
is a flat connection on $D^*$, then the conjugacy class of
$\phi_r(2\pi)$ is independent of $r$. In this case, denote the
conjugacy class of $\phi_r(2\pi)$ by $h_\a$. $h_\a$ is called the
{\it holonomy} of $\a$.

\begin{pro}\label{pro0}
Let $u=(u_1,u_2,\cdots,u_N)$ be a solution of (\ref{1.1})  in
$\R^n\backslash\{0\}$ with
\[u_i(x)=\mu_i\log |x| +O(1), \quad \text{  near } 0,\]
where $\mu_i>-2$ for $i\in I=\{1,2,\cdots, N\}$. Then the  flat
connection $\a$ defined by (\ref{eq3.5a}) has holonomy
\begin{equation}\label{holo}
h_\a= \left(\begin{matrix}e^{2\pi i \beta_0} & 0&\cdots &0\\
&  e^{2\pi i\beta_1} &0& \\
\cdots &\cdots & \cdots & \cdots \\
0 & 0& 0& e^{2\pi i \beta_N }\\
\end{matrix}\right),\end{equation}
where $\beta_0,\beta_1,\cdots,\beta_N$ are determined by
\begin{equation}\label{pro3.1.1} \beta_i-\beta_{i-1}=-\frac 12 \mu_i \quad  (i\in I)\quad
\text { and } \quad \sum_{j=0}^N\beta_j=0.\end{equation}
\end{pro}

For the proof, see \cite{JW2}.


Now by using Proposition  \ref{pro0}, we can generalize a result
in \cite{JW2}. Given an $N$-tuple $(\mu_1, \mu_2,\cdots,\mu_N)$
with $\mu_i>-2$ for any $i$ and let $u=(u_1,u_2,\cdots,u_N)$ be a
$C^2$ solution of the following system
\begin{equation}\label{eq2}
- \D u_i = \sum_{j=1}^N a_{ij} |x|^{\mu_j} e^{u_j}, \quad \hbox{
in } \R^2
\end{equation}
with
\begin{equation}\label{finite}
\int_{\R^2} |x|^{\mu_i} e^{u_i} <\infty.\end{equation}

Let
 \[m_i=\frac 1{2\pi}\int_{\R^2}|x|^{\mu_i}e^{u_i} \quad \hbox{ and } \quad
 \gamma_i=\sum a_{ij} m_j.\]
 The potential analysis gives that
 \begin{equation}\label{1}
 \begin{array}{llll}
 \ds\vs u_i&=&\ds -\gamma_i\log|x|+a_i+O(|x|^{-1}), & \quad \hbox{ near }\infty, \hbox{ and }\\
 \n u_i&=&\ds  -\gamma_i \frac{x}{|x|^2}+O(|x|^{-1})& \quad \hbox{ near }\infty,
 \end{array}\end{equation}
for some constants $a_i$ and
\begin{equation}\label{1_1}
\gamma_i-\mu_i>2,\end{equation} for all $i$. For the convenience
of the reader, we will give the proofs of (\ref{1} and (\ref{1_1})
in the appendix of the paper.

\begin{pro}\label{pro1} If $\mu_i\le 0$ for any $i$, then we have
\[\gamma_i=2(2+\mu_i),\quad \hbox{ for any } i.\]
\end{pro}

\medskip

\noindent{\it Proof.} Let $\tilde u=(\tilde u_1, \tilde
u_2,\cdots, \tilde u_N)$ with $\tilde u_i=u_i+\mu_i \log|z|$. It
is clear that $\tilde u$ satisfies (\ref{1.1}) in
$\R^2\backslash\{0\}$. Let $\a$ be the flat connection
(\ref{eq3.5a}) obtained from the solution $\tilde u$. We have
\[ \begin{array}{ll}
\tilde u_i(x)=\ds\vs \mu_i \log|x| +O(1), & \quad \hbox{ near } 0 \\
\tilde u_i(x)=\ds-(\gamma_i-\mu_i) \log|x| +O(1), & \quad \hbox{
near } \infty.
\end{array}\]
Let $\b_i$ be determined by $\{ \mu_i, i\in I\}$
 using (\ref{pro3.1.1}). By Proposition \ref{pro0}, the holonomy of $\a$
 at $0$ is $h_\a(0)$ given by (\ref{holo}).
 Now we compute the  the holonomy  of $\a$ at $\infty$, $h_\a(\infty)$.
To compute it, we use the Kelvin transformation to consider
\[ v_i(x)= \tilde u_i(\frac{\bar x}{|x|^2})-4\log|x|, \quad i\in I.\]
Clearly, $v=(v_1,v_2,\cdots,v_N)$ satisfies (\ref{1.1}) on
$\R^2/\{0\}$ and
\[ v_i=(\gamma_i-\mu_i-4)\log|x|+ O(1),\quad \hbox{ near } 0. \]
Let $\a'$ be its corresponding flat connection. It is obvious that
$h_\a(\infty)=h_{\a'}(0)$, which can be obtained again by
(\ref{holo}) by replacing $\mu_i$ by $\{\gamma_i-\mu_i-4\}$, i.e,
$h_{\a'}(0)$ is a matrix of form (\ref{holo}) with $\b'_i$ decided
by
\[\b_i'-\b_{i-1}'=-\frac 12 (\gamma_i-\mu_i-4)\quad \hbox{ and } \sum\b'_i=0.\]
Now the key fact is $h_\a(0)=h_\a(\infty)$, which implies that
\[\b_i-\b_i'= 1 \,\,\hbox{mod}\, {\mathbb Z}, \hbox{ for } i=0,1,2,\cdots, N.\]
Hence, $\gamma_i-2\mu_i=2n_i$ for some $n_i\in {\mathbb Z}$. By
(2.11), we have $n_i>1$ for any $i$.
 We now need a global Pohozaev identity, which is given in Proposition
 \ref{pro2} below. By (\ref{poho}), we have
\[\sum a^{ij}\gamma_i(2n_i-4)=0.\]
Since $a^{ij}>0$ for  any $i,j$ (see \cite{JW2}) and $n_i>1$,
 the previous formula implies that
$n_i=2$ for any $i$. This proves the Proposition. \qed

\begin{remark} To remove the condition $\mu_i\le 0$ for any $i$, one
 may have to relate the solutions of the Toda system to
 other geometric objects. For instance, a generalization of projective connections
 considered in \cite{Tro} might be a good candidate.
 \end{remark}
 \begin{pro}\label{pro2} $\gamma_i$ satisfies
 \begin{equation}\label{poho}\sum a^{ij}\gamma_i(\gamma_j-2(2+\mu_j))=0,
 \end{equation}
 where $(a^{ij})$ is the inverse matrix of the Cartan matrix $(a_{ij})$.
 \end{pro}
 \begin{proof} First for a given $R>0$, multiplying (\ref{1.1})
 by $x\cdot \n u_k$ and integrating over $B_R$, we have a Pohozaev identity
 \begin{equation}\label{2}\begin{array}{r}
 \ds\vs \frac 12\int_{B_R}\n u_k \n u_i +
 \frac 12 \int_{B_R}x\cdot \n (\n_ju_k) \n_ju_i\\
 \ds = \frac 12 \int_{\p B_R} R \frac{\p u_k}{\p n}\frac{\p u_i}{\p n}+\sum a_{ij}\int_{B_R}|x|^{\mu_j}e^{u_j} x\cdot \n u_k
 \end{array}\end{equation} From (\ref{2}), we have

 \begin{equation}\label{3}
 \begin{array}{c}
 \ds\vs \frac 12\sum a^{ki}\int_{B_R}\n u_k \n u_i +\frac 12
 \sum a^{ki}\int_{B_R}x\cdot \n (\n_ju_k) \n_ju_i
 \\
\ds  =-\frac 12 \sum a^{ki}\int_{\p B_R} R \frac{\p u_k}{\p
n}\frac{\p u_i}{\p n} +\sum \int_{B_R}|x|^{\mu_j}e^{u_j} x\cdot n
u_j.
 \end{array}\end{equation}
 It is easy to compute that
 \begin{equation}\label{4}
 \begin{array}{rcl} \ds \frac 12
\ds\vs  \sum a^{ki}\int_{B_R}x\cdot \n (\n_ju_k) \n_ju_i & =& \ds
\frac 14
\int_{B_R}x\cdot \n (\n_ju_k \n_ju_i) \\
 &=& \ds\vs  -\frac 12 \int_{B_R}\sum a^{ki}\n u_k \n u_i\\&&+\ds
 \frac 14 \sum a^{ik}\int_{\p B_R} R \n u_k\n u_i.
 \end{array}\end{equation}
 Inserting (\ref{4}) into (\ref{3}), we get
 \[
 \begin{array}{r}\ds\vs \sum \int_{B_R} (2+\mu_j)|x|^{\mu_j}e^{u_j}-\sum \int _{\p B_R}R|x|^{\mu_j}e^{u_j}
 \\ \ds =\frac 14 \sum a^{ik}(2\int_{\p B_R}\frac{ \p u_k}{\p n}
 \frac{ \p u_i}{\p n}-\int_{\p B_R} \n u_k\n u_i)
 \end{array}\]
Now using (\ref{1}) in the above formula and taking the limit
$R\to \infty$, we have proved the proposition .
 \end{proof}

 A geometric proof of the Proposition can be obtained
 by using  another geometric object, projective connections
 with singularities,
 which was used by Troyanov  \cite{Tro}
 to classify conical metrics of constant curvature.
 One can check that a quadratic differential  $\eta=fdz^2$, where
  $f:\R^2\slash\{0\}\to \C$ is given by
 \[f=\sum_{j,k=1}^N a^{jk}\{(\tilde u_{k})_{zz}-\frac 12 (\tilde u_j)_z\cdot
  (\tilde u_k)_z\}\]
 is a projective connection with regular
 singularities at $0$  and  $\infty$. By counting and comparing
 the weights of the
 singularities, one can get the global Pohozaev identity.

\section{A precise bubbling behavior}

In this section, we prove Theorem \ref{thm1.3}. Let
$u^k=(u_1^k, u_2^k)$ be a solution of
\begin{equation}
\label{eq001} \left\{ \begin{array}{rcl}
\ds -\D u_1^k & =&\vs \ds 2 h^k_1e^{u_1^k} -  h^k_2e^{u_2^k}\\
\ds-\D u_2^k& = &\ds2  h^k_2e^{u_2^k} -  h^k_1e^{u_1^k}\\
\end{array}\right. \ \ \mbox{in}\ \ B_2\end{equation}
satisfying conditions given in Theorem \ref{thm1.3}.
Set
\[\s_i=
 \lim_{k\to +\i} \int_{B_1} e^{u_i^k}.\]
 From Proposition \ref{pro}, we have
$\s_1=\s_2= 8\pi$.
%

In order to show the idea of our proof,  we first prove Theorem
\ref{thm1.3} in the case that $h^k_1=h^k_2=1$. Then we will point out which
steps should be modified for the general case.

\medskip

\noindent{\it Proof of Theorem \ref{thm1.3} in the case that $h^k_1=h^k_2=1$.}
 We divide the proof into several steps.
\medskip

\noindent{\bf Step 1.} From the above discussion, (\ref{1.9}) is valid in $B_{\ve_k R}$
for any fixed large number $R>0$. Hence, we only need to prove that (\ref{1.9})
is valid in $(\ve_kR,r_0)$ for some $r_0>0$. By assumption (4), we have
$$|u_2^k(x_1^k)-\l^k|<c,$$
for some constant $c>0$ independent of $k$. In the sequel, $c$
will
denote  a positive constant independent of $k$, which may vary
from line to line.
\medskip

\noindent
{\bf Step 2.}\quad By the Green representation, we have for $|x|<r_0$,
\begin{equation}{\label{eq:1-2}}\begin{array}{rcl}
 u_1^k(x) - \bar m_i^k
&=& \ds\vs \int_{B_{2r_0}} G(x,y) \left(2e^{u_1^k(y)} -
e^{u_2^k(y)}\right) dy +O(1)\\
&= &\ds
 \frac{1}{2\pi} \int_{B_{r_0}(x_\ve)} \log \frac{1}{|x-y|}
 \left(2e^{u_1^k(y)} - e^{u_2^k(y)}\right) dy +O(1),
\end{array}\end{equation}
where $\bar m_i^k=\inf_{\partial B_2}u^k_i(x)$.
In this paper, $O(1)$ denotes  a term bounded by a constant independent of $k$.
Hence, at $x_k$ we have
\begin{equation} {\label{eq:1-3}}
\l^k-\bar m_i^k = \frac{1}{2\pi}
\int_{B_{r_0}(x_\ve)} \log \frac{1}{|x_\ve-y|} \left( 2 e^{u_1^k(y)} -
e^{u_2^k(y)}\right) dy +O(1).
\end{equation}
(\ref{eq:1-2}) and (\ref{eq:1-3}) imply
\begin{eqnarray}{\label{eq:1-4}}
v_1^k(x) & = & u_1^k (x_k +\ve_k x) - u_1^k (x_k)\\ \nonumber
& = & \frac{1}{2\pi} \int_{B_{r_0}\ve_k^{-1}} \log \frac{|x|}{|y-x|}
\left(2e^{v_1^k(y)} - e^{v_2^k(y)}\right) dy +O(1).
\end{eqnarray}
Similarly, we have
$$u_2^k (x_k+\ve_kx)-u_2^k(x_k) = \frac{1}{2\pi}
\int_{B_{r_0}\ve_k^{-1}} \log \frac{|y|}{|y-x|}
\Biggl( 2e^{v_2^k(y)}-e^{v_1^k(y)}\Biggl) dy +O(1).$$
Since $|u_2^k(x_k)-\l^k|\leq c$, we have

\begin{equation}{\label{eq:1-5}}
v_2^k(x)=\frac{1}{2\pi} \int_{B_{r_0}\ve_k^{-1}} \log
\frac{|y|}{|y-x|} \left(2e^{v_2^k(y)} - e^{v_1^k(y)}\right) dy
+O(1).
\end{equation}

Set $w^k_1=\frac 13 (2v^k_1+v^k_2)$ and $w^k_2=\frac 13 (v^k_1+2v^k_2)$.
It is clear that
\[\begin{array}{rcl}
w^k_1&=&\vs\ds \frac{1}{2\pi} \int_{B_{r_0}\ve_k^{-1}} \log
\frac{|y|}{|y-x|} e^{v_1^k(y)} dy +O(1)\\
w^k_2&=&\ds \frac{1}{2\pi} \int_{B_{r_0}\ve_k^{-1}} \log
\frac{|y|}{|y-x|} e^{v_2^k(y)} dy +O(1).\end{array}\] A standard
potential analysis (see for instance \cite{CLT}) gives
\begin{equation}\label{eqa1}
\begin{array}{rcl} \ds \left|w^k_i+\frac
{\sigma^k_i}{2\pi}\log |x|\right| \le \d \log |x|+O(1),
\end{array}\end{equation}
for a small $\d>0$ and $|x|$ large enough, where $\s_i^k$ is the local mass
defined by
$$\s_i^k=\int_{B_{r_0}} e^{u_i^k}.$$ Thus, we have
\begin{equation}\label{eq-new1}
\begin{array}{rcl}
\ds \left|v^k_i-m^k_i\log |x|\right| \le 3\d \log |x|+O(1),
\end{array}\end{equation}
where $$m_1^k=\frac{2\s_1^k-\s_2^k}{2\pi} \quad \hbox{ and }\quad
m_2^k=\frac{2\s_2^k-\s_1^k}{2\pi}.
$$
As in \cite{CLT}, we choose a small constant $\d_3>0$ such that for
$\log (\frac{1}{\ve_k})\leq |x| \leq \frac{1}{\ve_k}$,
\begin{equation}\label{eq-new2}
\left|\tilde m^k_i(x)+m^k_i\right|=O(1)(\log \ve^{-1}_k)^{-1},
\end{equation}
where
\[\tilde m^k_i(x)=\int _{|y|\le\d_3|x|} e^{v_i^k(y)}dy.\]
This can be done, since from (\ref{eq-new1}) we have
\begin{equation}\label{eqa2}
\int_{|y|\ge\d_3\log \ve_k^{-1}} e^{v_i^k(y)}dy\le
 c\int_{|y|\ge\d_3\log \ve_k^{-1}} |y|^{-m^k_i+3\d}dy=
 O(1)(\log \ve_k^{-1})^{-3/2},\end{equation}
 by noting that $m^k_i\to 4$ as $k\to \infty$.
It is also easy to check that
\[\int_{|y|\ge\d_3|x|}\log \frac{|y|}{|x-y|} e^{v_i^k(y)}dy\le O(1)
(\log \ve_k^{-1})^{-1},\]
for $|x|\ge \log \ve_k^{-1}$.
Therefore, we have
\[\begin{array}{rcl}
w_i^k&=&\ds\vs\frac{1}{2\pi} \int_{|y|\le \d_3|x|} \log
\frac{|y|}{|y-x|} e^{v_1^2(y)} dy +O(1)
\\
&=&\ds\vs\frac{1}{2\pi} \int_{|y|\le \d_3|x|} \log
\frac{|y|}{|y-x|} e^{v_1^2(y)} dy +O(1)\\
&=&\ds\vs\tilde m^k_i\log|x|+O(1),
\end{array}\]
for $\log \ve^{-1}_k\leq |x| \leq \frac{1}{\ve_k}$.
Hence, we have
\begin{equation}{\label{eq:1-6}}
 \left| v_i^k(x) + m^k_i \log |x| \right|<c,
\end{equation}
for $\log \ve^{-1}_k\leq |x| \leq \frac{1}{\ve_k}$.
Similarly, we can show that
\begin{equation}{\label{eq:1-7}}
  \left| \nabla v_i^k(x) + m^k_i \frac{x}{|x|^2}\right|<c\\
\end{equation}
holds for $\log \ve^{-1}_k\leq |x| \leq \frac{1}{\ve_k}$.

By
(\ref{eq:1-6})-(\ref{eq:1-7}), we have
\begin{eqnarray}{\label{eq:1-7-0}}
\begin{array}{rcl}
\ds \vs u_i^k(x) &=& \ds m_i^k \log \frac{1}{|x|} + (2-m_i^k)\log \frac{1}{\ve_k}+O(1)\\
\ds \nabla u_i^k(x) &=&\ds  -m^k_i\frac{x}{|x|^2} +O(1)\end{array}
\end{eqnarray}
in $B_{r_0}\backslash B_{{\ve_k}\log \ve^{-1}_k}$.
\medskip

\medskip

\noindent
{\bf Step 3.} We want to show that
$$m_i^k=4\pi +O(1) (\log \frac{1}{\ve_k})^{-2}$$
for large $k$. Here we will use the geometric properties of the Toda system,
described in Section 3.

Let $\tilde w_0, \tilde w_1$ and $\tilde w_2$ defined by
\begin{equation}\label{eq100.1}
 2(\tilde w_i-\tilde w_{i-1}) = v^k_i \quad (i=1,2) \hbox{ and } \sum_{i=0}^2
 \tilde w_i=0.
\end{equation}
Let $w_1=\tilde w_1 -\log 4$ and $w_2=\tilde w_2 -2\log 4$.
Since $u_1,u_2$ satisfy the Toda system, we already know that
\begin{equation}{\label{eq:1-8}}
d\a+\frac 12 [\a, \a] =0,
\end{equation}
where $\a$ is defined by
$$\a := {\cal U}dz+{\cal V}d\bar z.$$
Here
$${\cal U} = \left(
\begin{array}{ccc}
(w_0)_z & e^{w_1-w_0} & 0 \\
0 & (w_1)_z & e^{w_2-w_1} \\
0 &0 & (w_2)_z \end{array}\right)$$
and
$${\cal V}= -\left(\begin{array}{ccc}
(w_0)_{\bar z}&0&0\\
e^{w_1-w_0} & (w_1)_{\bar z}&0\\
0 & e^{w_2-w_1}& (w_2)_{\bar z}\\
\end{array}\right).$$
It is clear that $\a$ is a one-form valued in $su(3)$. Equation (\ref{eq:1-8})
means that the connection $d+\a$ is flat.
Now we calculate the holonomy of $d+\a$ along
$$\{ |x|=(\log \frac{1}{\ve_k})^s\}$$
where $s>0$ will be fixed later. Decompose $\a$ by
$$\a=\a_r  dr +\a_\t d\t.$$
By definition, the holonomy $h_\a$ of the  connection $\a$ is
the conjugacy class of $g_r(2\pi)$, where $g_r(\theta)$ is the unique solution
of
\begin{equation}\label{eq101}
\begin{array}{rcl}
 \ds\vs \frac {d g_r}{d\theta}+\a_\theta g_r& =& 0,\\
\ds  g_r(0) & =&0.
\end{array}\end{equation}
 {From}  (\ref{eq:1-7-0}), 
it is easy to check that
\begin{equation}\label{eq120}
\begin{array}{rcl}
\a_\theta&=& \vs\ds \left( \begin{array}{ccc}
-i \b_0^k  &  b_1 ie^{i\theta} & 0\\
b_1  ie^{-i\theta}
& -i\b_1^k&  b_2 ie^{-i\theta}\\
0 & b_2ie^{i\theta} & -i\b_2^k\\
\end{array}
\right)\\
&+& \ds O(1)\epsilon_k\log\frac 1 {\epsilon_k} I\end{array}\end{equation}
at $|x|=(\log\frac{1}{\ve_k})^s$,
where \[b_{i}=O(1)
(\log \frac{1}{\ve_k})^{(1-\frac{m_i^k}{2})s}.\]
Since $m_i>2$ for all $i$, we choose $s>0$ such that
$s(\min_i m_i^k-2)>4$.

Now we can compute the holonomy $h_\a$ of $\a$ along $\{|x|=
 (\log \frac{1}{\ve_k})^s\}$
\begin{equation}\label{eq121}
\left( \begin{array}{ccc}
e^{-2\pi i{\b_0^k}}\\
& e^{-2\pi i{\b_1^k}}\\
& 0  & e^{-2\pi i  \b_2^k}
\end{array}\right)+O(1) (\log \frac 1{\epsilon_k})^{-2} I.\end{equation}
Since the holonomy $h_\a$ must be trivial, we have
$$\b_i^k=n_i+O(1) (\log \frac 1 \ve)^{-2}$$
for some integer $n_i$, which implies
\begin{equation}\label{eqa3}
m_i^k=2n'_i+O(1) (\log \frac{1}{\ve_k})^{-2}\end{equation} for
some integer $n'_i$. From Lemma \ref{lem3.1} we know $m^k_i\to 4$
as $k\to \infty$. Hence, we have
$$m_i^k=4+O(1)(\log \frac{1}{\ve_k})^{-2},$$
which,
together with (\ref{eq:1-6}), implies
\begin{equation}\label{eq-new5}\begin{array}{rcl}
v_i^k & = & \ds\vs-m_i^k \log |x| +O(1)\\
& = & \ds\vs-4\log |x| +O(1) (\log \frac{1}{\ve_k})^{-2} \log |x|+O(1)\\
& = & \ds -4 \log |x| +O(1),
\end{array}\end{equation}
for $\log\ve^{-1}_k\le |x| \le \frac 1{\ve_k}$.

\medskip

\noindent{\bf Step 4.} We have proved that (\ref{eq-new5}) holds
in $\log\ve^{-1}_k\le |x| \le \frac 1{\ve_k}$.
By Step 1, we have, for any
large $R>0$, $v^k_i$ converges to $v^0_i$ uniformly for $|x|\le R$.
By a classification result for entire solutions to the Toda system \cite{JW2},
we have
\[|v^0_i+4\log|x|| \le c\]
for $|x|\ge R$, where $c>0$ is a constant independent of $R$. Therefore,
\begin{equation}\label{eqb1}
|v^k_i+4\log|x|| \le 2c\end{equation}
for $|x|=R$ and large $k$. Choose $R$ large such that
\[e^{v^k_i} \le |x|^{-3}\]
for $|x|\ge R$ and $i=1,2$, and define
\[w_\pm (x)=-4\log|x|\pm(c_1-c_1|x|^{-\frac 12}).\]
It is clear that
\[\D w_\pm=\mp \frac 14 c_1|x|^{-\frac 52}.\]
 By choosing $R$ and $c_1$ large, we have, from the maximum principle,
 that
 \[w_-(x)\le v^k_i(x)\le w_+(x)\]
 for $R\le |x|\le \ve_k^{-1}$. Now we complete the proof of Theorem \ref{thm1.3}
 in the case $h^k_1=h^k_2=1$.
 \qed

\bigskip

 \noindent{\it Proof of
Theorem \ref{thm1.3} in the general case.} Assume without loss of
generality that $h_i^k(x)=h_i(x)$ and $h_i(0)=1$ for $i=1,2$. It
is clear that we only need to prove Step 3, i.e.,
\begin{equation}\label{eqn1}
m^k_i=4\pi+O(1)(\log \frac 1\ve)^{-2},\end{equation}
where
\[ m^k_1=\frac{2\s_1^k-\s_2^k}{2\pi} \quad\hbox{ and }\quad
m^k_2=\frac{2\s_2^k-\s_1^k}{2\pi}.\]
When $h_i$ is not constant, (\ref{eq001}) is not an integrable
system. However, the same idea as in Step 3 still works.

Define $\tilde w_i$ and $w_i$ as above.
We define a connection $\a=-({\cal U}dz+{\cal V}d\bar z)$ on the trivial bundle,
where
$${\cal U} = \left(
\begin{array}{ccc}
(w_0)_z & f^k_1e^{w_1-w_0} & 0 \\
0 & (w_1)_z & f^k_2e^{w_2-w_1} \\
0 &0 & (w_2)_z \end{array}\right)$$
and
$${\cal V}= -\left(\begin{array}{ccc}
(w_0)_{\bar z}&0&0\\
f^k_1e^{w_1-w_0} & (w_1)_{\bar z}&0\\
0 & f^k_2e^{w_2-w_1}& (w_2)_{\bar z}\\
\end{array}\right).$$
Here $f_i^k=h_i^{1/2}(\ve^kx+x^k)$.
Note that now $v^k=(v^k_1,v^k_2)$ satisfies in
$\O_k:=\{ x\in \R^2 \mid \ve_k x+x^k \in B_2\}$
\begin{equation}\label{eq102}
 \left\{ \begin{array}{l}
\ds -\D v_1^k = 2 h_1(\ve_k x+x^k)e^{v_1^k} - h_2(\ve_k x+x^k) e^{v_2^k},\\
\ds -\D v_1^k = 2 h_2(\ve_k x+x^k) e^{v_2^k} - h_1(\ve_k x+x^k)e^{v_1^k},\\
\end{array}\right.\end{equation}
 Now the connection $ d+\a$ is not flat. But it
satisfies
\[F=d\a+\frac 12[\a,\a]=\left(\begin{array}{ccc}
0&(f^k_1)_{\bar z} e^{w_1-w_0}&0\\
-(f^k_1)_ze^{w_1-w_0} & 0&(f^k_2)_{\bar z}e^{w_2-w_1} \\
0 & -(f^k_2)_z e^{w_2-w_1}& 0\\
\end{array}\right) dz\wedge d \bar z.\]
Here $F=F_{r\t}dr\wedge d\t$ is the curvature of the connection $d+\a$.

Let $g_r(\t)$ be the solution of (\ref{eq101}) for this connection $d+\a$.
We denote the conjugacy class of $g_r(2\pi)$ by $h_\a^r$. Now $h_\a^r$
may depend on $r$. Computing as in Step 3, we have (\ref{eq120}) and
\begin{equation}\label{eq130}
h_\a^r=\left( \begin{array}{ccc}
e^{-2\pi i{\b_0^k}} & 0&0\\
0& e^{-2\pi i{\b_1^k}}&0\\
0& 0  & e^{-2\pi i  \b_2^k}
\end{array}\right)+O(1) (\log \frac 1{\epsilon_k})^{-2} I,\end{equation}
at $r=(\log \frac 1{\ve_k})^s$.

Now we consider gauge transformations. Let $\phi:
\R^2\to SU(N+1)$. From the connection
$\a$, we can get a new connection by
\[ \tilde \a= \phi^{-1}d\phi+\phi^{-1}\a \phi.\]
For the new connection, we compute $h^r_{\tilde \a}$.
\begin{lemma} [\cite{SS}] \label{lem10}We have $h^r_{\tilde \a}=h^r_\a$.\end{lemma}
Using this Lemma, we can compute $h^r_\a$ by using a
suitable equivalent connection. This new connection is chosen as follows.
By solving in $\O_\ve\backslash B_R$
\[\begin{array}{rcl}
\ds\vs \frac {d\phi}{dr}+\a_r\phi &=&0\\
\ds \phi(0, \theta)&=& I, \end{array}\]
we obtain a new connection,
$\tilde \a=\phi^{-1}d\phi+\phi^{-1}\a \phi.$ Decompose the
new connection as above $\tilde \a=\tilde a_\t d\t+\tilde \a_r dr.$
Clearly, we have $\tilde \a_r=0$. Hence,
we have
\begin{equation}\label{eqc2}
\frac{d \tilde \a_\t}{dr}=\tilde F_{r\t},\end{equation}
where $\tilde F=\tilde F_{r\t}dr \wedge d\t$ is the curvature of the new connection. It is well-known
that $\tilde F=\phi^{-1} F\phi$.

Now we estimate $\tilde F_{r\theta}$.
We claim that
\begin{equation}\label{eqc1}
e^{v^k_1} \le C|x|^{-3}, \quad \hbox{ for } |x|\ge R,\end{equation}
for some fixed $R>0$
The claim will be proved at the end of the proof.
  {From} (\ref{eqc1}), it is clear that
 \[|(f^k_1)_ze^{w_1-w_0}|\le C\ve^k e^{\frac 12 v^k_1} \le C \ve^k |x|^{-3/2},\]
 for $|x|\le R$.
  We have similar estimates for
the other entries.
 Therefore, we have
 \[|\tilde F_{r\t}| =|F_{r\t}|\le C \ve^k |x|^{-1/2}
 \quad \hbox{ for } |x|\ge R.\]
 On $B_R$ we have $|\tilde F_{r\t}| =|F_{r\t}|\le C\ve^k$. Therefore, from
 (\ref{eqc2}) we have
 \[|\tilde \a_\t(r,\t)|\le C\ve(\log\frac{1}{\ve^k})^sI, \quad \hbox{ for }
  r=(\log \frac 1{\ve^k})^s.\]
  Now it is clear to see that
  \[h^r_{\tilde \a}=O(1)\ve(\log\frac{1}{\ve^k})^sI,\]
  which, together with (\ref{eq130}) and Lemma \ref{lem10}, implies
  \[\b^k_i=n_i+O(1)(\log\frac 1{\ve^k})^{-2},\]
  and hence
  \[m^k_i=4+O(1)(\log\frac 1{\ve^k})^{-2}.\]
 This finishes the proof of Step 3 and the Theorem  for the general case.

 Now it remains to check  Claim (\ref{eqc1}). The idea  is similar to Step 4.
 We define for $|x|\ge R$ a function $w$
 \[w =-3\log |x|+(c_1-c_1|x|^{-\frac12}).\]
 We have $\D w=-\frac1 4 c_1|x|^{-\frac 52}$. Since $m^k_i\to 4$ as
 $k\to \infty$, we may assume $m^k_i>3 $ for large $k$. Applying the
 maximum principle, we get the claim. \qed

\section{Existence: A critical case}

 Let $\S$ be a Riemann surface with
Gaussian curvature $K$. In this section and the next section, we will
 study the existence of solutions of
\begin{equation}
\label{eq6.1}\begin{array}{rcll} -\D u_i &=&\ds\vs \sum^N_{j=1}
\rho_ja_{ij}\left(\frac{h_je^{u_j}} {\int_\S h_je^{u_j}
}-1\right),&\quad \hbox{ in }
 \Omega\quad 1\le i\le N,\\
\end{array}
\end{equation}
for the  coefficient matrix $A=(a_{ij})_{n\times n}$, the  Cartan
matrix of $SU(N+1)$ and $\rho=(\rho_1,\rho_2,\cdots, \rho_N)$ with
$\rho_i>0$ ($i=1,2,\cdots, N$) given  constants. System
(\ref{eq6.1}) is the Euler-Lagrange system of the
functional
\begin{equation}\label{eq6.2}
J_\rho(u)=\frac12\sum_{i,j=1}^N\int a^{ij}\n u_i\n
u_j+\sum_{j=1}^N\int \rho_ju_j-\sum_{j=1}^N \rho_j
\log\int_{\S}h_je^{u_j},
\end{equation}
in $H:=(H^1(\S))^N.$ When $\rho_i<4\pi$, it
is a subcritical case  with respect to the Moser-Trudinger inequality
established in \cite{JW1}. By the Moser-Trudinger inequality, the 
functional $J_\rho$
is coercive and hence has  a minimizer.
When $\rho_i \le 4\pi$ and  some of
 the
$\rho_i$'s are equal to $4\pi$,  the functional $J_\rho$ still
has a lower bound, but it does not satisfy the coercive condition.
It is a critical case and the existence problem becomes more subtle. Even
in the case that $N=1$, the existence problem is a difficult problem, see
\cite{DJLW, NT, ChL}. Here we will use results in \cite{ChL} to obtain
a sufficient condition under which the functional $J_\rho$ achieves its minimum.
For simplicity we only consider the case $N=2$.

\begin{proposition} \label{pro6.1} For a fixed $\rho_2\in (0, 4\pi)$, define
 \[\cal X=\{\hbox{Solutions of (\ref{1.6}) for } \rho_1\in (0,4\pi]
\hbox{ and } \rho_2\}.\] If (\ref{1.5}) holds, then $\cal X$
is a compact set.
\end{proposition}
\begin{proof} If $\cal X$ is not compact, then we may assume that there is
a sequence $\{\rho_k\}$ with $\rho_k<4\pi$ and $\rho_k \to 4\pi$
as $k\to \infty$, and a sequence of solutions $u^k=(u^k_1, u^k_2)$
satisfying (\ref{1.6}) with $\rho_1=\rho_k$ and $\rho_2$ such
that
\[\max\{\max u^k_1(x), \max u^k_2(x)\} \to \infty\]
as $k\to \infty$. Since $\rho_2<4\pi$, by Proposition 2.5, we
know that $\max u^k_1 \to \infty$ as $k\to \infty$ and $u_2^k$ is
uniformly bounded from above. 
Let $w^k$ be the
unique solution of
\begin{equation}\label{eq-1}
\begin{array}{rcl}
-\D w^k &= & \ds\vs \rho_2\left(\frac{h_2e^{u^k_2}}{\int_\Sigma
h_2e^{u^k_2}}-1\right)=:f_k
\end{array}\end{equation}
with $\int w^k=0$.
Since $u^k_2$ is uniformly bounded from above, $w^k$ is uniformly
bounded in $H^{2,q}$ for any $q<\infty$. Hence $w^k$ converges to
$w_0$ in $C^{1,\a}$ for some $\a>0$.

Let $ v^k_1=u^k_1+w^k$ and $v^k_2=u^k_2-2w^k$. It is clear that
$v^k_1$ and $v^k_2$ satisfy
\begin{equation}\label{eq6.4.a}
-\D v^k_1= 2\rho_k\left(\frac{\tilde h^k_1 e^{v^k_1}}{\int_\Sigma
\tilde h^k_1 e^{v^k_1}}-1\right)
\end{equation} and
\[-\D v^k_2=\ds - \rho_k \left(\frac{\tilde h^k_1 e^{v^k_1}}
{\int_\Sigma \tilde h^k_1 e^{v^k_1}}-1
\right).\]with $\tilde h^k_1=h_1e^{-w^k}$.
It is clear that $v^k_1+2v^k_2=0$.
Set $\tilde v^k_1=v^k_1-\log\int_\Sigma \tilde h_1^ke^{v^k_1}$.
It is clear that $\tilde v^k_1$ satisfies
\begin{equation}\label{eq6.4}
-\D \tilde v^k_1= 2\rho_k\left(\tilde h^k_1 e^{\tilde v^k_1}-1\right).
\end{equation}
In order to apply Theorem 1.1 in \cite{ChL} to our situation, we need to show
that $\tilde h^k_1=h_1e^{-w^k}$ converges to $h_1e^{-w^0}$ in $C^{2,\a}$.
 Let $p$ be the unique blow-up
point of $u^k_1$. As in \cite{JW1}, one can show that $ v^k_1\to
8\pi G(p,\cdot)$ and $ v^k_2 =-v^k_1/2\to -4\pi G(p,\cdot)$ in
$H^{1,q}(\Sigma)$ ($q<2$) and in
$C^2_{loc}(\Sigma\backslash\{p\})$, where $G(x, \cdot)$ is
defined by
\[-\D G(x,\cdot) =\d_x-1\]
and $\int G=0$. We apply a result of \cite{Li} (see also
\cite{CLT}) to (\ref{eq6.4}) and obtain that
\begin{equation}\label{eq-2}
|\tilde v^k_1(x)-U_k(x) |<c,
\quad \hbox{ for } |x-p|<r_0\end{equation}
for some constant $c>0$ and a small constant $r_0>0$. Here
\[U_k(x)=\log\frac{e^{\l_k}}{(1+\frac
{e^{\l_k}h_1(p)e^{-w_0(p)}}{8}|x|^2)^2}\]
and $\l_k=\max \tilde v^k_1$.
Furthermore, one can prove that
\begin{equation}\label{eq-2.0}
|\l_k-\log \int_\Sigma \tilde h_1^k e^{v_1^k}|<c,\end{equation}
and
\begin{equation}\label{eq-2.1}
|\n \tilde v_1^k (x)-\n U_k(x)|<c\end{equation}
for some constant $c>0$ and $|x-p|\leq r_0$. See [4].
{From} (\ref{eq-2}), (\ref{eq-2.0}), (\ref{eq-2.1}) and $u^k_2=-\frac 12 v^k_1+2w^k$, we have
\[|\n e^{u^k_2}|\le c
e^{-\frac12 \l_k}|\n e^{-\frac12 \tilde v^k_1}|
 \le c_1(1+ e^{-\frac12 \l_k}|\n e^{-\frac12 U_k}|)\le c_2.\]
Therefore, the elliptic estimates implies that
$w^k$ converges in $C^{2,\a}$ to $w^0$
satisfying
\begin{equation}\label{eq44}
-\D w^0=\ds \rho_2\left(\frac{h_2e^{-4\pi G}e^{2w^0}} {\int_\Sigma
h_2e^{-4\pi G}e^{2w^0}} -1\right).
\end{equation}

One can check by applying the Pohozaev identity that the blow-up
point $p$ satisfies
\begin{equation}\label{eq5}
\n \log h_1(p)+\n \tilde G(p)-\n w^0(p)=0,\end{equation}  where
$\tilde G$ is the regular part of the Green function $G$. Note
that $\D w^0(p)=\rho_2$. Since $w^k$ converges to $w^0$ strongly in $C^{2,\a}$,
we can apply Theorem 1.1 in \cite{ChL2}
to (\ref{eq6.4}) and obtain that
\begin{equation}\label{eq6.5}
2\rho_k-8\pi =2h_1^{-1}(p)e^{w^0(p)}(\D \log
h_1(p)+8\pi-\rho_2-2K(p))\l_ke^{-\l_k}+o(e^{\l_k}),\end{equation}
where $\l_k=\max \tilde v^k_1$, which tends to $\infty$. Since $2\rho_k<8\pi$,
(\ref{eq6.5}) leads to a contradiction to (\ref{1.5}). Hence $\cal X$ is compact.
\end{proof}

\noindent{\it Proof of Theorem \ref{thm1.1}.} Choose a sequence
$\{\rho^k\}$ with $\rho^k<4\pi$ and $\rho^k\to 4\pi$. From above,
we know that for $\rho^k=(\rho^k, \rho_2)$ with $\rho_2<4\pi$ the
functional $J_\rho^k$ has a minimizer $u^k$. By Proposition
\ref{pro6.1}, we have that $u^k$ converges to $u^0$, which is
clearly a minimizer of $J_\rho$ with $\rho=(4\pi, \rho_2)$.\qed

\begin{theorem}\label{thm6} Let $\S$ be a Riemann surface with
Gaussian curvature $K$.
Suppose that
\begin{equation}\label{1.5a}
\min\{\D \log h_1(x), \D \log h_2(x)\}+4\pi -2K(x)>0.\end{equation}
 Then $J_\rho$ has a minimizer
$u=(u_1,u_2)$ satisfying (\ref{1.6}) for $\rho=(\rho_1,\rho_2)=(4\pi, 4\pi)$.
\end{theorem}

\begin{proof} Let $\rho^k_2<4\pi$ be a sequence tending to $4\pi$.
Applying Theorem \ref{thm1.1}, we know that $J_\rho$ has a
minimizer $u^k=(u^k_1, u^k_2)$ satisfying (\ref{1.6}) for
$\rho=(4\pi, \rho^k_2)$. To prove the Theorem, we only need to show that
$u^k$ converges strongly to $u^0$. Assume by contradiction that
$u^k$ blows up, namely
\[\max\{\max_{\Sigma} u^k_1,\max_{\Sigma} u^k_2\}\to \infty.\]
Let
\[S_i=\{x\in\Sigma\,|\, \exists \hbox{ a sequence }\{x_k\} \hbox{ with }
\lim_{k\to \infty}x_k=x \quad \& \quad
\lim_{k\to \infty}u^k_i(x_k)=\infty\}.\]
By Proposition 2.5, we know that $S_1 \cap S_2=\emptyset$. Hence we have two possibilities:
\begin{itemize}
\item [(1)] $S_1=\{p\}$ and $S_2=\emptyset$
or $S_2=\{p\}$ and $S_1=\emptyset$.
\item [(2)] $S_1=\{p\}$ and $S_2=\{q\}$ with $p\neq q$.
\end{itemize}
For case (1), using the argument given in the proof of Proposition
\ref{pro6.1}, we can compute $8\pi-\rho^k_2$ in terms of the
blow-up speed as in (\ref{eq6.5}) to get a contradiction
to condition (\ref{1.5a}). For case (2), we consider
a disk $D$ on $\Sigma$ with $p\in D$ and $q\not\in \overline D$. $u^k=(u^k_1, u^k_2)$
satisfies on $D$
\begin{equation}\label{eq6.3-1}
\begin{array}{rcl}
-\D u_1&=&\ds\vs 8\pi\left(\frac {h_1e^{u_1}}{\int_\S
h_1e^{u_1}}-1\right)
-\rho_2\left(\frac {h_2e^{u_2}}{\int_\O h_2e^{u_2}}-1\right)\\
-\D u_2&=&\ds 2\rho_2\left(\frac {h_2e^{u_2}}{\int_\S
h_2e^{u_2}}-1\right) -4\pi\left(\frac {h_1e^{u_1}}{\int_\S
h_1e^{u_1}}-1\right).
\end{array}\end{equation}
On the boundary $\partial D$, we know that there are two sequences $\{a^k_1\}$
and $\{a^k_2\}$ with $a^k_i\to -\infty$ for $i=1,2$ such that
$u^k_i-a^k_i$ is bounded. Define
\[-\D w^k=\rho^k_2
\left(\frac{h_2 e^{u^k_2}}{\int_\Sigma h_2 e^{u^k_2}}-1\right)\]
with the Dirichlet boundary condition $w^k(x)=0$ for $x\in \partial D$. Now
consider $v^k_1=u^k_1+w^k-a^k_1$ which satisfies
\[-\D v^k_1=8\pi \left(\frac{\tilde h^k_1 e^{v^k_1}}
{\int \tilde h^k_1 e^{v^k_1}}-1\right),\]
where
\[\tilde h^k_1=e^{a^k_1}h_1^k e^{-w^k}.\]
Note that ${v^k_1}_{|\partial D}$ is bounded.
Define $v^k_0$ by
\[\left\{\begin{array}{rcl}
\vs-\D v_0^k&=& 8 \pi \quad\hbox{ in } D\\
 v^k_0&=& -v^k_1 \quad\hbox{ on }\partial  D
\end{array}\right.\]
and set $\tilde v^k =v^k_1+v^k_0$. It is clear that $\tilde v^k_1$ satisfies
\[
\left\{
\begin{array}{rcl} -\D \tilde v^k_1&=&\vs\ds
8\pi \frac{h_k e^{\tilde v^k_1}}
{\int_\Sigma h_k e^{\tilde v^k_1}}\quad\hbox{ in } D\\
\tilde v^k_1 &=& 0 \quad \hbox{ on }\partial D
\end{array}\right. \]
with $h_k=\tilde h^k_1e^{-v^k_0}$. It is also clear that
\[\frac{\int_D h_k e^{\tilde v^k_1}} {\int_\Sigma h_k e^{\tilde v^k_1}}
\to 1,\]
as $k\to \infty$. As in the proof of Proposition \ref{pro6.1}, we can
show that $h_k$ converges strongly in $C^{2,\a}$. Hence, again we can use
the result in \cite{ChL}  for the Dirichlet boundary problem
to get that the divergence of $\tilde v^k_1$
implies that
\begin{equation}\label{66}
\lim_{k\to\infty}(\D \log h_1-\D v^k_0-\D w^k-2K)=0\end{equation}
at the blow-up point $p$. Since we have $-\D w^k \to -4\pi$
and $-\D v^k_0=8\pi$,  from (\ref{66}) we have
\[ \D \log h_1+4\pi-2K=0,\]
which contradicts  (\ref{1.5a}).
\end{proof}

\

A direct consequence of Theorem \ref{thm1.1} is
\begin{corollary} Let $\Sigma$ be a closed surface of genus
greater than $0$. Suppose that $\Sigma$ has constant curvature and $h_1$
 and $h_2$ are constants.
Then equation (\ref{1.6}) has a solution for any $\rho=(\rho_1,\rho_2)$
with $\rho_i \le 4\pi$.
\end{corollary}


\section{Existence: A supercritical case}
In this section we consider the case that $\rho_i\in (0, 8\pi)$,
but $\rho_i\neq
4\pi$ for any $i$. This is a supercritical case, in the sense that
the functional (\ref{1.4}) has no lower bound.
The result obtained in this section is not related to the curvature
of the underlying surface. Therefore, in this section we consider
the existence in a bounded domain $\O\subset \R^2$. A similar result
as obtained in Theorem \ref{thm2} holds for a closed surface with genus
greater that $0$, which is Theorem \ref{thm1.4}.

\begin{theorem}\label{thm2}
Let $\O$ be a smooth bounded domain and
 $\rho_i\in (0,4\pi)\cup (4\pi,8\pi)$ $(i=1,2)$ two constants. Suppose that the boundary
 $\partial \O$ of $\O$ has at least two components.
For any positive
function $h_i$ $(i=1,2)$,
there exists a solution $u=(u_1,u_2)$ satisfying
\begin{equation}\label{eq5.3}
\begin{array}{rcl}
\D u_1&=&\ds\vs 2\rho_1\frac {h_1e^{u_1}}{\int_\O h_1e^{u_1}}
-\rho_2\frac {h_2e^{u_2}}{\int_\O h_2e^{u_2}}\\
\D u_2&=&\ds 2\rho_2\frac {h_2e^{u_2}}{\int_\O h_2e^{u_2}}
-\rho_1\frac {h_1e^{u_1}}{\int_\O h_1e^{u_1}},
\end{array}\end{equation}
and $u_i=0$ on $\partial \O$ $(i=1,2)$.
\end{theorem}

When $N=1$, the result was proven in \cite{DJLW2}, see also
\cite{ST}. Recently,  Lucia and Nolasco \cite{LN} obtained a
non-trivial solution of (\ref{eq5.3}) under some conditions
when $h_i$ is a constant for $i=1,2$.

Our method of the proof closely follows \cite{DJLW2}, which, in
turn, is motivated by \cite{St} (see also \cite{St-book}). Hence,
here we only sketch the main ideas of the proof. When one of
the
$\rho_i$ is $4\pi$, the existence becomes more subtle. We will
consider this case elsewhere.

\medskip

\noindent{\it Proof of Theorem \ref{thm2}.} When $\rho_i<4\pi$
$(i=1,2)$, the Theorem is true for any bounded domain. We first
consider the case $\rho_i\in (4\pi, 8\pi)$ $(i=1,2)$. We divide
the proof into several steps.

\medskip

\noindent{\bf Step 1.} We first define the center of mass of a function
$v\in H^1_0(\O)$  by
\[{m_c}(v)=\frac{\int_{\O}xe^v}{\int_{\O}e^v}.\]
Assume, for the simplicity of the notation, that
$\partial \O=\partial_+ \O\cup\partial_- \O$ has only two disjoint
components. Define a family of functions
\begin{equation}
\label{eq5.4}
\Gamma:(-\infty, +\infty)\to
H^1_0(\O)\times H^1_0(\O)\end{equation}
with $\Gamma(t)=(\gamma_1(t),\gamma_2(t))$ satisfying
\begin{equation}
\label{eq5.5} J_\rho(\Gamma(t)) \to \pm \infty \quad \hbox{ as }
t\to - \infty
\end{equation}
and
\begin{equation}
\label{eq5.6} m_c(\gamma_i) \to \partial_{\pm} \O\quad \hbox{ as } t\to\pm \infty.
\end{equation}
The existence of such a family is guaranteed by $\rho_i>4\pi$
$(i=1,2)$. Denote the set of all such families by $\cal X$ and
define a minimax value
\begin{equation}
\label{eq5.7} \a_\rho:=\inf _{\Gamma\in {\cal
X}}\sup_{t}J_\rho(\Gamma(t)).\end{equation}

\medskip

\noindent{\bf Step 2.} The minimax value $\a_\rho>-\infty$.

\medskip

The proof of the step follows from an improved   Moser-Trudinger
inequality under a condition introduced by Aubin.

\begin{lemma}\label{lem1}
Let $S_1$ and $S_2$ be two subsets of $\bar{\O}$ satisfying
$dist(S_1,S_2)\ge\delta_0>0$ and $\gamma_0\in (0,1/2)$.
For any $\e>0$, there exists a constant $c=c(\e,\delta_0,\gamma_0)>0$
such that
$$J_{(8\pi-\ve, 8\pi-\ve)}(u) \ge -c$$
holds for all $u=(u_1,u_2)\in H^{1}_0(\O)\times H^{1}_0(\O)$ satisfying
\begin{equation}\label{eq5.8}
\frac{\int_{S_1}e^{u_i}}{\int_{\O}e^{u_i}}\ge \gamma_0
\qquad \hbox{and}\qquad
\frac{\int_{S_2}e^{u_i}}{\int_{\O}e^{u_i}}\ge \gamma_0,
\end{equation}
for $i=1,2$.
\end{lemma}

Let $\gamma_0\subset \O$ be a closed curve in $\O$ enclosing the inner boundary of
$\O$. Each curve starting from
$\partial_-\O$ and ending at $\partial_+\O$ intersects with $\gamma_0$.
By Lemma \ref{lem1}, we can show that
\[J_\rho(u)>-c,\]
for any $u\in H^1_0(\O)\times H^1_0(\O)$ with center of mass
 $m_c(u)\in \gamma_0$.
See the argument in \cite{DJLW2}. Hence we prove the step.

\medskip

\noindent{\bf Step 3.} $\a_\rho:(4\pi, 8\pi)\times (4\pi, 8\pi)$
is non-increasing in the following sense: if
$\rho=(\rho_1,\rho_2)$ and $\rho'=(\rho_1',\rho'_2)$ with $\rho_1
\le \rho_1'$ and $\rho_2=\rho_2'$, then
\[\frac{\a_{\rho}}{\rho_1} \ge\frac{\a_{\rho'}}{\rho'_1}.\]

\medskip

This is easy to check.
\medskip

\noindent{\bf Step 4.} Now fix $\rho_2$ and define
\[\L_1=\L_1(\rho_2)=\{\rho_1\in (4\pi, 8\pi)\,|\,\frac{\a_{\rho_1,\rho_2}}{\rho_1}
\hbox{ is differentiable at }\rho_1\}.\] It is clear that $\L_1$
is dense, i.e., $\bar \L_1=[4\pi, 8\pi]$. Following a method given
by Struwe \cite{St} and \cite{St-book}, see also \cite{DJLW}, we
can prove that for fixed $\rho_2$ and $\rho_1\in \L_1$,
$\a_{\rho_1,\rho_2}$ is achieved by $u=(u_1,u_2)$, which is a
solution of (\ref{eq5.3}).
\medskip

\noindent{\bf Step 5.} Now for any $\rho=(\rho_1, \rho_2) \in
(4\pi,8\pi)\times  (4\pi,8\pi)$, we have a sequence
$\rho^k=(\rho^k_1, \rho_2)$ with $\rho^k_1\in\L_1$ and
$u^k=(u^k_1,u^k_2)$ satisfying (\ref{eq5.3}) with $\rho^k$.

Since $\rho_2\neq 4\pi m$, by Proposition \ref{pro} we see that $u^k$
converges  to $u^0=(u^0_1, u^0_2)$, which is a solution of
(\ref{eq5.3}). This finishes the proof for the case $\rho_i\in
(4\pi,8\pi)$ $(i=1,2)$.

\medskip

When $\rho_1< 4\pi$ and $\rho_2\in (4\pi,8\pi)$, $\Gamma$ defined
by (\ref{eq5.4})-(\ref{eq5.6}) does not exist. We need to modify
the definition. Consider $\Gamma=(\gamma_1,\gamma_2)$ with
(\ref{eq5.4}) and
\begin{equation}
m_c(\gamma_2) \to \partial_{\pm}\O \quad \hbox{ as } t\to\pm \infty.
\end{equation}
Then the same argument finishes the proof of the Theorem.\qed

Theorem \ref{thm1.4} can be proven in a similar way.

Furthermore, in this case
 we can compute the topological degree of the solution space
as in \cite{ChL2}. First we define the Leray-Schauder degree for
system (\ref{1.6}). For the Leray-Schauder degree for
the Liouville type equation, see \cite{Li} and \cite{ChL2}.
Consider $\rho=(\rho_1, \rho_2)$ with $\rho_i\neq 0 \, \hbox{mod}\, 4\pi$.
Define an operator
\[T_\rho:=\D^{-1} \left(\begin{matrix} 2 &-1\\
-1 & 2 \end{matrix}\right)\left(\begin{matrix} \rho_1 & 0\\
0 & \rho_2 \end{matrix}\right)\ds \left(\begin{matrix} \frac{h_1 e^{u_1}}
{\int h_1 e^{u_1}}-1  \\ \frac{h_2 e^{u_2}}
{\int h_2 e^{u_2}}-1
 \end{matrix}\right).\]
 $T_\rho$ acts on $H^1_0\times H^1_0$, where $H^1_0=\{u\in H^1\,|\, \int u=0\}.$
Set ${\cal X}_\rho$ the solution space of (\ref{1.6}).
 {From} above, we know that ${\cal X}_\rho$ is compact. Hence we can
 define the Leray-Schauder degree for (\ref{1.6})
 \[d_\rho:=\deg(I+T_\rho, B_R,0),\]
 where $I$ is the identity and $B_R:=\{u=(u_1,u_2)
  \in H^1_0\times H^1_0\,|\, \|u_1\|+\|u_2\|\le R\}$ for a large $R>0$.
It is clear that $d_\rho$ is well-defined. The homotopy
invariance of the Leray-Schauder degree implies that
 $d_\rho$ is independent of $h_1$ and $h_2$. Furthermore Theorem \ref{thm1.2}
 implies that
 $d_{\rho}=d_{\tilde \rho}$ if there are two integers $m_1$ and $m_2$ such that
 $\rho_i \in (4 \pi (m_i-1), 4\pi m_i)$ and
 $\tilde \rho_i \in (4 \pi (m_i-1), 4\pi m_i)$.

 \begin{theorem}\label{thm10} Let $\Sigma$ be a closed surface
 of genus $g$. Then we have for $\rho=(\rho_1, \rho_2)$ with
 $\rho_1 \in(0,4\pi)$ and $\rho_2\in (4\pi, 8\pi)$
 \[d_\rho=2g-1.\]
 \end{theorem}

  The proof follows from the argument given in the proof of Theorem
   \ref{thm1.1} and a result given in \cite{ChL2}. In a forthcoming paper,
   we will compute the topological degree for the Toda system for the
   general case.

\section{Appendix. Behavior of singularities of  solutions with finite energy}
In this Appendix, we consider the asymptotic behavior of
singularities of  solutions to the Toda system.

Let $u=(u_1,u_2,\cdots, u_N)$  be a solution of
\begin{equation}\label{eq7.1}
-\D u_i =\sum_{i=1}^Na_{ij} e^{u_j}, \quad \hbox{ for } i=1,2,\cdots, N
\end{equation}
in a punctured disk $D^*=D\backslash\{0\}$ with
\begin{equation}\label{eq7.2}
\int_D e^u_i <\infty.\end{equation}
By the potential analysis, one can check that there is
a constant $\gamma_i$
\begin{equation}\label{eq7.3}
\lim_{|x|\to 0}\frac{u_i}{-\log|x|}=\gamma_i,\end{equation} for
each $i$. (\ref{eq7.2})  and (\ref{eq7.3}) imply that $\gamma_i\le
2$ for any $i$. Furthermore, one can check that
\begin{equation}\label{eq7.4}
\lim_{|x|\to 0}\frac{u_i+\gamma_i\log|x|}{\log|x|}=0.
\end{equation}

\begin{proposition}\label{pro7.1}
For any $i$, we have
\[  \gamma_i<2.\]
\end{proposition}

\noindent{\it Proof.}
 Suppose that there is $j \in I$ such that $\gamma_{j}=2$.
We  claim that $\gamma_{j-1}=2$ or $\gamma_{j+1}=2$.
Assume by contradiction that $\gamma_{j-1}<2$ and $\gamma_{j+1}<2$.
We assume that neither $j =1$
nor  $j=N$. (For the cases $j=1$ and $j=N$, the argument is the same.)
   Consider in $D^*$
\[-\D u_j=-e^{u_{j-1}}+2e^{u_j}-e^{u_{j+1}}=:f.\]
Since $\gamma_j=2$, $\gamma_{j-1}<2$ and $\gamma_{j+1}<2$, we know
that $f(x)>0$ in a small punctured disk $D^*_\d=D_\d\backslash\{0\}$.
Set
\[v(x)=-\frac 1{2\pi}\int_{D_\d}\log|x-y| f(y)\]
and $w=u_j-v$. It is clear that $\D v=-f$ and $\D w=0$.
One can check that
\[ \lim_{|x|\to 0} \frac {v}{-\log|x|}=0,\]
which implies
\[\lim_{|x|\to 0}\frac{w(x)}{-\log|x|}
=\lim_{|x|\to 0}\frac{u_j-v}{-\log|x|}=2.\]
Since $w$ is harmonic in $D_\d^*$, we have
$w=-2\log|x|+w_0$ with a smooth harmonic function $w_0$ in $D_\d$.
By definition, we know $v>0$. Thus, we have
\[\int_{D_\d}e^{u_j}=\int_{D_\d}e^{w+v}\ge
\int_{D_\d}\frac 1{|x|^2}e^{w_0}=\infty,\]
a contradiction. A similar argument implies that
if $\gamma_{j}=\gamma_{j+1}=\cdots =\gamma_{j+k}=2$, then
either $\gamma_{j-1}=2$ or $\gamma_{j+k+1}=2$, by using the equation
for $\frac 1k(u_{j}+u_{j+1}+\cdots +u_{j+k})$. Hence we can show $\gamma_i=2$
for all $i$.
Now we consider $\tilde u=\frac 1 N(\sum_{i=1}^N u_i)$.
It is clear that $\tilde u$ satisfies
\[-\D \tilde w=\frac{1}{N}(e^{u_1}+e^{u_2})>0\]
with $\int_{D}e^{\tilde w}<\infty$ and $\lim_{|x|\to 0}\frac{\tilde w}{-\log|x|}=2.$
The same argument as above gets a contradiction. Hence $\gamma_i<2$ for any
$i$.
\qed

 \
 
\noindent{\it Proof of (\ref{1}) and (\ref{1_1}).}
Now we apply the Proposition to solutions of (\ref{eq2}) with condition
(\ref{finite}). Let $u=(u_1,u_2,\cdots, u_N)$ be such a solution.
The potential analysis gives us
\begin{equation}\label{eq8.1}
\lim_{|x|\to \infty} \frac{u_i}{\log|x|}=-\gamma_i=-\frac 1 {2\pi}
\sum_{j=1}^Na_{ij}\int_{\R^2}|x|^{\mu_j}e^{u_j}.\end{equation}
We consider $\tilde u=(\tilde u_1, \tilde u_2,\cdots, \tilde u_N)$
with $\int_D e^{\tilde u_i}<\infty$ for any $i$ and
\[\tilde u_i=u_i(\frac x{|x|})-(4+\mu_i)\log|x|.\]
It is easy to check that $\tilde u$ satisfies (\ref{eq7.1}) in
$D^*$ with
\[\lim_{|x|\to 0}\frac{\tilde u_i(x)}{-\log|x|}=4+\mu_i-\gamma_i.\]
Applying the Proposition, we have
\[4+\mu_i-\gamma_i<2,\]
for any $i$. This is (\ref{1_1}).
 {From} the above inequalities and the potential analysis, we have
\[\tilde u_i=-(4+\mu_i-\gamma_i)\log|x|+O(1),\]
from which we have the first formula in (\ref{1}). The second
formula follows from the potential analysis. \qed

\begin{lemma}\label{lem-A} Let $u\in C^2$ be a solution of
\begin{equation}\label{eqA1}
-\D u= 2|x|^{2\a} |x-q|^2e^u\end{equation}
 with $\a\ge 0$ and
 \begin{equation}\label{eqA2}
 \gamma:=\frac 1{2\pi}\int_{\R^2}|x|^{2\a}|x-q|^2e^u<
 \i.\end{equation}
 then $\gamma > 2+\a.$\end{lemma}
 \begin{proof} We show this lemma following closely \cite{ww}.
 Set
 \[w=\frac
 1{2\pi}\int_{\R^2}(\log|x-y|-\log(|y|+1))2|y|^{2\a}|y-q|^2dy.\]
 It is easy to check
 by the potential analysis that
 $u=w+c$ for some constant $c$ and
 \[w \ge -2\gamma \log|x|-c_1,\]
for some constant $c_1>0$. This, together with the finiteness
 condition (\ref{eqA2}), implies that $\gamma >2+\a$.\end{proof}
 \begin{remark} In fact, all solutions of (\ref{eqA1}) and (\ref{eqA2})
 can be
 classified. \end{remark}

\end{document}